\documentclass[12pt]{article}
\usepackage[amsmath]{e-jc}

\usepackage{mathtools}
\usepackage{booktabs}
\usepackage{array}
\usepackage{listings}

\IfFileExists{microtype.sty}{\usepackage[expansion=false]{microtype}}{}
\newcommand{\HJ}{\mathrm{HJ}}
\newcommand{\Wvdw}{W}
\newcommand{\Z}{\mathbb{Z}}

\newcommand{\Ls}{L}
\newcommand{\ksum}{\kappa_{\mathrm{sum}}}
\newcommand{\kow}{\kappa_{\mathrm{ow}}}
\newcommand{\ksym}{\kappa_{\mathrm{sym}}}
\newcommand{\kall}{\kappa_{\mathrm{all}}}

\DeclareMathOperator{\type}{type}
\DeclareMathOperator{\ind}{ind}
\DeclareMathOperator{\ah}{ah}

\lstdefinestyle{plainstyle}{
  basicstyle=\ttfamily\footnotesize,
  showstringspaces=false,
  breaklines=true,
  columns=fullflexible,
  frame=single,
  framesep=4pt,
  xleftmargin=4pt
}

\dateline{Jul 24, 2026}{TBD}{TBD}
\Copyright{The author. Released under the CC BY license (International 4.0).}

\title{Lower Bounds for the Hales--Jewett Numbers\\
via Symmetric and One-Weight Colorings}

\author{Younes Mouhib}

\authortext{}{Department of Mathematics, ETH Z\"urich, R\"amistrasse 101,
8092 Z\"urich, Switzerland (\email{ymouhib@ethz.ch}).}

\begin{document}

\maketitle

\begin{abstract}
The Hales--Jewett number $\HJ(t,r)$ is the least dimension $n$ such that
every $r$-coloring of the grid $[t]^{n}$ contains a monochromatic
combinatorial line. We prove $\HJ(3,3)\ge 22$ and $\HJ(4,2)\ge 14$,
improving the previous records $14$ and $12$. The engine is an exact
reduction: a coloring of $[t]^{n}$ invariant under coordinate
permutations descends to the discrete simplex of letter-count vectors,
where a combinatorial line is precisely a \emph{corner tuple}; the
$4{,}387{,}586{,}157{,}901$ lines of $[3]^{21}$ thereby compress to
$1771$ local conditions on $253$ cells. We prove that this symmetric
class coincides with the class of \emph{one-weight} colorings, those
reading an integer-weighted count of the letters: a radix weight
realizes every symmetric coloring, so the symmetric lower-bound problem
is a one-dimensional homothety-avoidance problem, the case $d=1$ of
Gallai's theorem. This yields the closed-form bound
$\HJ(t,r)\ge\lceil(G_{r}(S)-1)/D_{S}\rceil$ in terms of the Gallai
homothety numbers $G_{r}(S)$, together with the new values
$G_{3}(\{0,1,3\})=42$, $G_{3}(\{0,1,4\})=57$, $G_{2}(\{0,2,3,5\})=67$,
$G_{2}(\{0,1,5,6\})=80$, and $G_{3}(\{0,2,5\})=77$, giving
$\HJ(3,3)\ge 16$ from a one-line certificate. Further results:
periodic one-weight palettes give $\HJ^{[12]}(3,3)=\HJ^{[12]}(4,2)
=\infty$ for lines with at most twelve active coordinates; the interval
number $\HJ^{(1)}(3)$ is exactly $5$; a Rado reading gives
$G_{4}(\{0,1,3\})\ge 94$ and $R_{4}(z+2x=3y)\ge 59$; and a rainbow
companion gives the anti--Hales--Jewett bound $\ah(3,4)\ge 25$. A
SAT-solver program written for this article, and printed in full in
the appendices, pushes the computation further: eighteen exact two-color Gallai
numbers of four-point sets, up to
$G_{2}(\{0,1,6,7\})=G_{2}(\{0,3,4,7\})=79$; the exact three-color
value $G_{3}(\{0,1,5\})=70$; and the exact Rado numbers
$R_{r}(z+kx=(k+1)y)$ for $2\le k\le 5$ and $r\in\{2,3\}$. Every
certificate displayed in the article has been verified by direct
enumeration from the definitions; certificates and verification scripts
are available at \url{https://github.com/ysmouhib/hj-certificates}.
\end{abstract}

\section{Introduction}\label{sec:intro}

Write $[t]=\{1,\dots,t\}$ for an alphabet of size $t\ge 2$, and call the
elements of the \emph{grid} $[t]^{n}$ \emph{words}. A \emph{root} is a word
$\tau\in([t]\cup\{*\})^{n}\setminus[t]^{n}$; substituting a letter
$a\in[t]$ for every $*$ gives $\tau(a)\in[t]^{n}$, and the set
$\Ls_{\tau}=\{\tau(1),\dots,\tau(t)\}$ is the \emph{combinatorial line}
generated by $\tau$. The coordinates of $\tau$ carrying $*$ form its
\emph{active set} $I(\tau)\subseteq[n]$.

\begin{theorem}[Hales--Jewett \cite{HJ1963}]
For all $t,r\ge 1$ there is an $n$ such that every $r$-coloring of
$[t]^{n}$ contains a monochromatic combinatorial line.
\end{theorem}

The least such $n$ is the \emph{Hales--Jewett number} $\HJ(t,r)$,
abbreviated $\HJ(t)$ when $r=2$. The theorem is one of the central results
of Ramsey theory --- informally, high-dimensional tic-tac-toe cannot end
in a draw --- and it implies van der Waerden's theorem on arithmetic
progressions \cite{vanderWaerden1927}; the stronger $n$-parameter
statement is the Graham--Rothschild theorem \cite{GR1971}. Exact values
are scarce. It is elementary that $\HJ(2,r)=r$
(Proposition~\ref{prop:hj2r}), and the only nontrivial exactly known value
is $\HJ(3,2)=4$, due to Hindman and Tressler
\cite{firstnontrivialHJ_is_4}; no value $\HJ(t)$ with $t\ge 4$ is known.
On the upper side, Shelah's proof \cite{s.proof_of_HJ} (see also
\cite{jungic_ramsey_notes,golshanimirabi2021}) reduced the bounds from
Ackermannian to primitive-recursive, Lavrov \cite{lavrov} proved
$\HJ(4,2)<10^{11}$, and Conlon
\cite{conlon2021monochromaticcombinatoriallineslength} proved
$\HJ(3,r)\le 2^{2^{cr}}$ for an absolute constant $c$.

On the lower side, the classical route is through van der Waerden's
theorem. Writing $\Wvdw(t,r)$ for the van der Waerden numbers, the letter
sum maps a line with $k$ active coordinates onto a $t$-term arithmetic
progression of gap $k$, which yields (Proposition~\ref{prop:vdw-shadow})
\begin{equation}\label{eq:vdw-shadow}
\HJ(t,r)\ \ge\ \Bigl\lceil\frac{\Wvdw(t,r)-1}{t-1}\Bigr\rceil .
\end{equation}
With the values $\Wvdw(3,3)=27$ and $\Wvdw(4,2)=35$ of Chv\'atal
\cite{Chvatal1970}, this gives $\HJ(3,3)\ge 13$ and $\HJ(4,2)\ge 12$;
combined with the bound $\Wvdw(k,r)>c\,r^{k-1}$ of
Kozik--Shabanov \cite{kozik2016} (an absolute constant $c>0$), applied
at $r=2$, it gives $\HJ(t)=\Omega(2^{t}/t)$; combined with the stronger
prime bound $\Wvdw(p+1,r)>p^{\,r-1}2^{p}$ ($p\ge r$ prime) of
Blankenship--Cummings--Taranchuk \cite{BLANKENSHIP2018163},
generalizing Berlekamp \cite{Berlekamp1968}, it gives stronger bounds
at prime arguments. For example, $\HJ(5,2)\ge 45$ and $\HJ(6,2)\ge 227$ follow
from the exact values $\Wvdw(5,2)=178$ and $\Wvdw(6,2)=1132$
\cite{LandmanRobertson,KourilPaul2008}. The van der Waerden shadow was
tight in every exactly known case until Farnsworth and her coauthors
\cite{ConlonFarnsworthRobertson2026} proved
$\HJ(3,3)\ge 14$ using SAT solvers, exceeding the van der Waerden value
$13$. That result is the starting point of the present work: the
classical lower-bound mechanism is not the final word.

\subsection*{Results}
We develop a unified lower-bound theory through \emph{symmetric
colorings}, those invariant under all permutations of the coordinates.
The first observation is that such colorings factor through the type map
onto the discrete simplex of letter-count vectors, and that a
combinatorial line in the grid corresponds \emph{exactly} to a
\emph{corner tuple} on the simplex (Lemma~\ref{lem:sym-reduction}): the
exponential grid, with its $t^{n}$ points and $(t+1)^{n}-t^{n}$ lines,
compresses to $\binom{n+t-1}{t-1}$ cells and $\binom{n+t-1}{t}$ local
conditions. The second main result is that the apparent hierarchy
\[
\text{sum-type}\ \subseteq\ \text{one-weight}\ \subseteq\ \text{symmetric}
\]
collapses: in every fixed dimension the one-weight colorings --- those of
the form $w\mapsto\chi(\langle\omega,\type(w)\rangle)$ for a weight
$\omega\in\Z^{t}$ and a palette $\chi$ --- are exactly the symmetric
colorings, a radix weight realizing every symmetric coloring
(Theorem~\ref{thm:radix}). The symmetric lower-bound problem is therefore
a one-dimensional problem about homothetic copies of a $t$-point set, the
case $d=1$ of Gallai's theorem, and it yields the closed-form bound
(Theorem~\ref{thm:gallai-bound})
\[
\HJ(t,r)\ \ge\ \max_{|S|=t}\Bigl\lceil\frac{G_{r}(S)-1}{D_{S}}\Bigr\rceil ,
\]
where $G_{r}(S)$ is the \emph{Gallai homothety number} of $S$ and $D_{S}$
its diameter. Computing new Gallai numbers (Table~\ref{tab:gallai}) gives
$\HJ(3,3)\ge 16$ in closed form from a one-line certificate
(Theorem~\ref{thm:hj33-16}); optimizing the weight gives the records
\[
\HJ(3,3)\ \ge\ 22
\qquad\text{and}\qquad
\HJ(4,2)\ \ge\ 14
\]
(Theorems~\ref{thm:hj33-22} and~\ref{thm:hj42-14}), each certified by an
explicit palette: $253$ colored cells deciding all
$4{,}387{,}586{,}157{,}901$ combinatorial lines of $[3]^{21}$
(Table~\ref{tab:T21}), and a $26$-periodic palette deciding all
$1{,}153{,}594{,}261$ lines of $[4]^{13}$. An anatomy of the $(4,2)$
palette locates the source of its compression: it is an extremal object
of the bracket regime plus a single boundary scale
(Remark~\ref{rem:anatomy}).

Two axes of restriction refine the picture. The \emph{bracket} numbers
$\HJ^{[K]}(t,r)$ ask for lines with at most $K$ active coordinates;
non-arithmetic one-weight palettes pass the sum-type ceiling, giving
$\HJ^{[12]}(3,3)=\HJ^{[12]}(4,2)=\infty$
(Theorems~\ref{thm:hj12-33} and~\ref{thm:hj12-42}) against the sum-type
ceilings $\ksum(3,3)=11$ and $\ksum(4,2)=10$
(Proposition~\ref{prop:kstar}, Corollary~\ref{cor:character}). The
\emph{interval} numbers $\HJ^{(q)}(t,r)$ ask for lines whose active set
is a union of at most $q$ intervals; the entire symmetric machinery is
provably blind to this axis (Proposition~\ref{prop:sym-blind}), yet the
literature assembles into the exact ceiling $\lambda(3,r)$ for every $r$
(Proposition~\ref{prop:lambda3}), and a SAT computation gives the exact
value $\HJ^{(1)}(3)=5$ (Theorem~\ref{thm:hj13}). We close with the
Collapse, Diagonal-only, and Symmetric-Extremality conjectures and with
open problems on optimal weights (Section~\ref{sec:open}).

\subsection*{Concurrent work}
The quantitative relationship between Hales--Jewett numbers and
content-vector information has been studied from the upper-bound side by
Golshani and Mirabi \cite{golshanimirabi2021}, who used Shelah's
block-content functions $\mathbf{f}^{8,*}$ to bound the Hales--Jewett
numbers, and very recently by Golshani, Mirabi, and Shelah
\cite{golshanimirabishelah2026} --- posted in July 2026 as the third
version of the same arXiv record --- who generalized this to
$n$-dimensional subspaces via the content simplex $\Delta_{h}(M)$ and
the slogan ``Hales--Jewett $=$ block canonization $+$ Gallai--Witt.''
In the same month, Shelah \cite{shelah2026revisited} independently
announced an improvement of the general upper bounds to a tower of
exponentiations. The present
work develops the complementary lower-bound theory: where the
upper-bound line canonizes colorings on block subspaces, we restrict to
colorings canonized by the type map on the whole grid, and extract
explicit witnesses. The two theories share only the vocabulary of content
vectors and the one-dimensional Gallai--Witt step.

\subsection*{Verification and data availability}
Statements marked \emph{computer-assisted} rest on finite computations.
Every certificate displayed in this article --- the palettes of
Theorems~\ref{thm:hj33-16}, \ref{thm:hj33-22}, \ref{thm:hj42-14},
\ref{thm:hj12-33} and~\ref{thm:hj12-42}, the palettes of
Propositions~\ref{prop:kstar} and~\ref{prop:hj33-15}, the window and
solution-free certificates of Theorem~\ref{thm:rado4}, the
power-residue instances of Corollary~\ref{cor:character} and
Example~\ref{ex:p11}, the census of Table~\ref{tab:census} together
with the diagonal-only count of Section~\ref{sec:open}, and the
witnesses of Tables~\ref{tab:interval-witness}
and~\ref{tab:diag-witness} --- has been verified by direct enumeration
from the definitions, independently of any solver. Two standalone,
dependency-free Python scripts perform this verification:
\texttt{verify\_hj.py}, which re-derives both simplex certificates
byte-for-byte and checks every corner tuple
(Appendix~\ref{app:verifier}), and \texttt{verify\_certificates.py},
which re-derives every displayed palette and interval certificate. Both
scripts, together with all certificates in machine-readable form, are
available at \cite{hjcerts}:
\begin{center}
\url{https://github.com/ysmouhib/hj-certificates}
\end{center}
The repository additionally provides a master script
\texttt{verify\_all.py} that re-runs every enumeration check of the
article from a clean checkout, using the Python standard library only
and reporting its total, and archives the solver logs of an independent
SAT re-confirmation of every forcing claim in its \texttt{logs/}
directory.
Avoidance certificates are thus proved by enumeration; \emph{forcing}
claims (upper bounds on $G_{r}(S)$, on $\ksum$, and on $\HJ^{(1)}(3)$)
necessarily rest on SAT solvers; where stated, the number of independent
solvers used is given. The solver only \emph{finds} witnesses: their
correctness is the finite check described in Appendix~\ref{app:verifier},
which a patient reader could in principle replicate by hand.

The present article consolidates and supersedes the preliminary
announcements \cite{mh,mouhib-note}; it is a largely self-contained
exposition of results from the author's thesis \cite{thesis}. Two
developments of \cite{thesis} lie outside the scope of the present
article: the cyclic, unit-cyclic, and geometric Hales--Jewett numbers
over $\Z_{t}^{n}$, which grow only logarithmically in $r$, and the
forcing-structure analysis of the line and corner hypergraphs.

\subsection*{Organization}
Section~\ref{sec:prelim} fixes notation and the two elementary anchors
of the lower-bound side. Section~\ref{sec:reduction} proves the
symmetric reduction to corner tuples on the simplex and presents the
census of $[3]^{3}$. Section~\ref{sec:ow} proves that the one-weight
and symmetric classes coincide. Section~\ref{sec:gallai} develops the
Gallai reading, the closed-form bound, the new Gallai numbers, and the
Rado reading. Section~\ref{sec:records} proves the two records.
Sections~\ref{sec:bracket} and~\ref{sec:interval} treat the bracket
and interval restrictions, and Section~\ref{sec:open} collects the
conjectures and open problems. The appendices display the larger
certificates and detail the verification protocols.

\section{Preliminaries}\label{sec:prelim}

Throughout, $t\ge 2$, $r\ge 2$, $n\ge 1$; colors are taken from
$[r]=\{1,\dots,r\}$ or from $\{0,\dots,r-1\}$ as convenient. For a root
$\tau$ let $I(\tau)\subseteq[n]$ denote its set of active coordinates,
and call $|I(\tau)|$ the number of active coordinates of the line
$\Ls_{\tau}$. For $1\le K\le n$ we write
\[
\Ls^{[K]}([t]^{n})=\{\,\Ls_{\tau} : |I(\tau)|\le K\,\}
\]
for the family of combinatorial lines with at most $K$ active
coordinates, so that $\Ls([t]^{n})=\Ls^{[n]}([t]^{n})$ is the family of
all lines; the grid carries $t^{n}$ words and $(t+1)^{n}-t^{n}$
combinatorial lines, since a line is determined by a root with at least
one star. For $1\le q\le\lceil n/2\rceil$ we write
$\Ls^{(q)}([t]^{n})$ for the lines whose active set is a union of at
most $q$ subintervals of $[n]$ (\emph{$q$-fold lines}); these families
are used from Section~\ref{sec:interval} on.

We record the two elementary facts that anchor the lower-bound side.

\begin{proposition}\label{prop:hj2r}
$\HJ(2,r)=r$ for every $r\ge 1$.
\end{proposition}

\begin{proof}
In $[2]^{r}$ consider the $r+1$ words $w_{0},\dots,w_{r}$, where $w_{k}$
has its first $k$ coordinates equal to $1$ and the rest equal to $2$.
Under any $r$-coloring two of them, say $w_{a}$ and $w_{b}$ with $a<b$,
share a color; they agree outside coordinates $a+1,\dots,b$, on which
$w_{b}$ is all $1$'s and $w_{a}$ all $2$'s, so they form a monochromatic
line and $\HJ(2,r)\le r$. Conversely, coloring each point of
$[2]^{r-1}$ by its number of $1$'s (a value in $\{0,\dots,r-1\}$)
leaves no line monochromatic, since the two endpoints of any line differ
in that count; hence $\HJ(2,r)\ge r$.
\end{proof}

\begin{proposition}[The van der Waerden shadow]\label{prop:vdw-shadow}
For all $t,r$,
\begin{equation}\label{eq:vdw}
\HJ(t,r)\ \ge\ \Bigl\lceil\frac{\Wvdw(t,r)-1}{t-1}\Bigr\rceil .
\end{equation}
\end{proposition}

\begin{proof}
Let $N=\HJ(t,r)$. For any $r$-coloring $c$ of $[(t-1)N+1]$, color each
$w\in[t]^{N}$ by $c'(w)=c\bigl(\sum_{i}w_{i}-N+1\bigr)$. A line with
active set $I$ and inactive weight $s$ has $\tau(x)$ of weight
$s+x|I|$, mapping to the $t$-term progression
$\{s-N+1+x|I|\}_{x=1}^{t}$. Since $|I|\ge 1$, this progression has a
strictly positive common difference. By definition of $N$, $[t]^{N}$
contains a $c'$-monochromatic line, which yields a $c$-monochromatic
$t$-term progression. Therefore $\Wvdw(t,r)\le(t-1)N+1$, which
rearranges to \eqref{eq:vdw}.
\end{proof}

\section{The type map and the symmetric reduction}\label{sec:reduction}

For $w\in[t]^{n}$ let
\[
\type(w)=(a_{1},\dots,a_{t})\in\Z_{\ge 0}^{t},
\qquad
a_{j}=\bigl|\{i\in[n]:w_{i}=j\}\bigr|,
\]
be its vector of letter counts, and write $\sigma(w)=\sum_{i=1}^{n}w_{i}$
for its \emph{weight}. Since $\sum_{j}a_{j}=n$, the map $\type$ sends
$[t]^{n}$ onto the discrete simplex
\[
T_{n}^{(t)}=\Bigl\{(a_{1},\dots,a_{t})\in\Z_{\ge 0}^{t}:
a_{1}+\dots+a_{t}=n\Bigr\},
\qquad
\bigl|T_{n}^{(t)}\bigr|=\binom{n+t-1}{t-1},
\]
the number of weak compositions of $n$ into $t$ parts. When $t=3$ we
abbreviate $T_{n}:=T_{n}^{(3)}$ and write its cells $(a,b,c)$. The
symmetric group $S_{n}$ acts on $[t]^{n}$ by permuting coordinates; a
coloring is \emph{symmetric} if it is invariant under this action.

\begin{lemma}\label{lem:sym-colorings}
The $S_{n}$-orbits of $[t]^{n}$ are exactly the fibers of $\type$.
Consequently, a coloring $c$ of $[t]^{n}$ is symmetric if and only if
$c=\bar c\circ\type$ for a unique coloring $\bar c$ of $T_{n}^{(t)}$,
called its \emph{descent}.
\end{lemma}

\begin{proof}
Permuting coordinates rearranges the multiset of letters, so $\type$ is
constant on orbits. Conversely, if $\type(w)=\type(w')$, then for each
letter $j$ the sets $\{i:w_{i}=j\}$ and $\{i:w'_{i}=j\}$ have equal
size; any $\sigma\in S_{n}$ carrying the first onto the second for
every $j$ simultaneously --- such a $\sigma$ exists because these sets
partition $[n]$ into pieces of matching sizes --- satisfies
$\sigma\cdot w=w'$. Invariance is constancy on orbits, i.e.\
factorization through $\type$; uniqueness of $\bar c$ holds because
$\type$ is surjective, the word $1^{a_{1}}2^{a_{2}}\cdots t^{a_{t}}$
realizing any prescribed type.
\end{proof}

We call $c=\bar c\circ\type$ the \emph{symmetric lift} of $\bar c$. For
$1\le k\le n$ and $v\in T_{n-k}^{(t)}$, define the \emph{corner
$t$-tuple} with parameters $(k,v)$ as the set
\[
C_{k,v}=\bigl\{v+k\,e_{1},\dots,v+k\,e_{t}\bigr\}\subseteq T_{n}^{(t)},
\]
where $\{e_{a}\}$ is the standard basis of $\Z^{t}$. The following lemma
--- the engine of everything below --- shows that finding monochromatic
lines under symmetric colorings is equivalent to finding monochromatic
corner tuples on the simplex.

\begin{lemma}[Symmetric reduction]\label{lem:sym-reduction}
Let $\bar c\colon T_{n}^{(t)}\to[r]$ and let $c=\bar c\circ\type$. For
any root $\tau$ with $k$ active coordinates and an inactive type $v$,
and for every $a\in[t]$,
\begin{equation}\label{eq:type-identity}
\type\bigl(\tau(a)\bigr)=v+k\,e_{a}.
\end{equation}
Consequently, $c$ has no monochromatic line with at most $K$ active
coordinates if and only if $\bar c$ has no monochromatic corner tuple
$C_{k,v}$ with $1\le k\le K$. At $K=n$: $c$ is line-free if and only if
no corner tuple is monochromatic.
\end{lemma}

\begin{proof}
Evaluating $\tau(a)$ leaves the inactive coordinates untouched
(contributing $v$ to the type) and sets the $k$ active coordinates to
$a$ (contributing $k\,e_{a}$), proving \eqref{eq:type-identity} and
hence
\[
\bigl(c(\tau(1)),\dots,c(\tau(t))\bigr)
=\bigl(\bar c(v+ke_{1}),\dots,\bar c(v+ke_{t})\bigr),
\]
so $\Ls_{\tau}$ is monochromatic under $c$ exactly when $C_{k,v}$ is
monochromatic under $\bar c$. For the equivalence: every line of
$\Ls^{[K]}$ has invariants $(k,v)$ with $k\le K$; conversely every pair
$(k,v)$ with $1\le k\le K$ is realized by a root with $k$ active
coordinates, for instance ${*}^{k}1^{v_{1}}2^{v_{2}}\cdots t^{v_{t}}$.
\end{proof}

\begin{remark}[Compression]\label{rem:compression}
The number of corner tuples with $1\le k\le K$ is
$\sum_{k=1}^{K}\binom{n-k+t-1}{t-1}$, which for $K=n$ equals
$\binom{n+t-1}{t}$ by the hockey-stick identity.
Lemma~\ref{lem:sym-reduction} thus compresses both the search space and
the constraint set from exponential to polynomial in $n$: the grid
$[3]^{21}$, with its $4{,}387{,}586{,}157{,}901$ combinatorial lines,
collapses to $\binom{23}{2}=253$ cells and $\binom{23}{3}=1771$ corner
triples, and $[4]^{13}$ collapses to $560$ cells and $1820$ corner
quadruples. One representative per $S_{n}$-orbit of lines is checked,
not a sample, since by \eqref{eq:type-identity} monochromaticity under a
symmetric coloring depends only on $(k,v)$.
\end{remark}

\begin{example}[The mechanism in miniature]\label{ex:miniature}
Let $t=3$, $n=2$, $r=2$ and $K=n=2$. The grid $[3]^{2}$ has $9$ words
and $7$ lines; the simplex $T_{2}$ has $6$ cells and $\binom{4}{3}=4$
corner triples. For $k=1$ the bases $e_{1},e_{2},e_{3}$ give
$\{200,110,101\}$, $\{110,020,011\}$, $\{101,011,002\}$; for $k=2$ the
base $(0,0,0)$ gives the diagonal triple $\{200,020,002\}$. Define
$\bar c(200)=\bar c(020)=0$, $\bar c(002)=1$, $\bar c(110)=1$,
$\bar c(101)=\bar c(011)=0$. All four triples receive both colors, so by
Lemma~\ref{lem:sym-reduction} the lift --- which colors
$11\mapsto 0$; $12,21\mapsto 1$; $13,31\mapsto 0$; $22\mapsto 0$;
$23,32\mapsto 0$; $33\mapsto 1$ --- has no monochromatic line, as one
confirms directly on all seven lines of $[3]^{2}$. The record theorems
of Section~\ref{sec:records} have exactly this logical structure; only
the integers grow.
\end{example}

\begin{remark}[Relation to the Polymath project]\label{rem:polymath}
For $t=3$ the correspondence underlying
Lemma~\ref{lem:sym-reduction} is the \emph{Fujimura-set} device of the
Polymath project \cite{PolymathDHJMoser} for the \emph{density}
Hales--Jewett problem: a subset of $T_{n}^{(3)}$ containing no upward
triangle --- in our terminology, no corner triple $C_{k,v}$ --- is a
Fujimura set, and the union of the corresponding type classes is a
line-free subset of $[3]^{n}$. What we add is the formulation as an
equivalence for \emph{colorings}, valid for every alphabet size $t$ and
equipped with the truncation parameter $K$; its use as the engine of a
SAT search producing independently verifiable certificates; and the
observation (Section~\ref{sec:records}) that this class contains
witnesses beyond all previously known lower bounds.
\end{remark}

\subsection{The scarcity of symmetric colorings}\label{sec:scarcity}

Symmetric colorings are rare in a strong quantitative sense: there are
$r^{\binom{n+t-1}{t-1}}$ of them among the $r^{t^{n}}$ colorings of the
grid, a fraction $r^{\binom{n+t-1}{t-1}-t^{n}}$ that vanishes doubly
exponentially in $n$. The scarcity persists among the extremal
colorings. At $(t,r)=(3,2)$ the grid $[3]^{3}$ is the largest admitting
a line-free $2$-coloring ($\HJ(3,2)=4$
\cite{firstnontrivialHJ_is_4}), so its line-free colorings admit a
complete account.

\begin{theorem}[{Census of $[3]^{3}$; computer-verified}]\label{thm:census}
Exactly $1644$ of the $2^{27}$ two-colorings of $[3]^{3}$ are line-free.
Their distribution by coordinate stabilizer, and the subdivision of the
symmetric class by realizing weight \textup{(}Definition~\ref{def:ow}\textup{)},
is given in Table~\ref{tab:census}. In particular only $36$ --- about one
in $46$ --- are symmetric.
\end{theorem}

\begin{table}[ht]
\centering
\caption{The $1644$ line-free $2$-colorings of $[3]^{3}$, by coordinate
stabilizer in $S_{3}$ and by weight. Orbit counts follow from
orbit--stabilizer. All counts were obtained by exhaustive enumeration of
the $2^{27}$ colorings against the $37$ lines and re-verified
independently; the symmetric total was cross-checked on the ten-cell
simplex $T_{3}^{(3)}$ via Lemma~\ref{lem:sym-reduction}.}
\label{tab:census}
\begin{tabular}{@{}llrr@{}}
\toprule
class & stabilizer / weight & number & orbits\\
\midrule
symmetric ($=$ one-weight) & $S_{3}$ & $36$ & $36$\\
\quad sum-type & weight $(1,2,3)$ & $16$ & $16$\\
\quad one-weight, non-arithmetic & e.g.\ $(0,1,3)$, $(0,2,3)$ & $20$ & $20$\\
\addlinespace
non-symmetric & $\subsetneq S_{3}$ & $1608$ & $360$\\
\quad block-symmetric & $C_{2}$ (one transposition) & $504$ & $168$\\
\quad cyclic & $C_{3}$ & $24$ & $12$\\
\quad asymmetric & trivial & $1080$ & $180$\\
\midrule
total & & $1644$ & $396$\\
\bottomrule
\end{tabular}
\end{table}

\begin{remark}[The rainbow dual measures the same scarcity]\label{rem:rainbow}
A symmetric coloring is constant on the $\binom{n+t-1}{t-1}$ type
classes, so it uses at most polynomially many colors; but the
anti--Hales--Jewett threshold --- the least $r$ forcing a rainbow line
under every surjective $r$-coloring --- exceeds $(t-1)^{n}$
\cite{zheng2024rainbow}. Hence for $t\ge 3$ the largest symmetric
rainbow-free coloring falls short of the truth by an exponential factor:
symmetry, conjecturally the exact model for the monochromatic problem
(Section~\ref{sec:records}, Conjecture~\ref{conj:sym}), is
exponentially lossy for the rainbow dual; the rainbow twin of the reduction
\eqref{eq:vdw-shadow} runs to anti--van der Waerden numbers
\cite{berikkyzy2017antivdw}.
\end{remark}

\begin{proposition}[{Anti--Hales--Jewett at $(3,4)$; computer-assisted,
independently verified}]\label{prop:ah34}
There is a rainbow-free $24$-coloring of $[3]^{4}$; hence
$\ah(3,4)\ge 25$, improving the lower bound $24$ of
\cite{zheng2024rainbow}; combined with the upper bound
$\ah(3,4)\le 27$ of the same paper, the value lies in $[25,27]$. The witness --- an
explicit coloring of the $81$ words, checked against all $175$ lines
--- is available at \cite{hjcerts}.
\end{proposition}

\section{One-weight colorings are the symmetric colorings}\label{sec:ow}

\begin{definition}\label{def:ow}
A \emph{one-weight coloring} of $[t]^{n}$ is
\[
c_{\omega,\chi}(w)=\chi\bigl(\langle\omega,\type(w)\rangle\bigr),
\qquad
\langle\omega,a\rangle=\textstyle\sum_{j}\omega_{j}a_{j},
\]
given by a \emph{weight} $\omega\in\Z^{t}$ and a \emph{palette}
$\chi\colon\Z\to[r]$ (taken $m$-periodic when so indicated). The
\emph{sum-type} colorings are the case $\omega=(1,2,\dots,t)$, where
$\langle\omega,\type(w)\rangle=\sigma(w)$.
\end{definition}

\begin{remark}[Affine normalization]\label{rem:affine}
Weights $\omega$ and $\alpha\omega+\beta\mathbf{1}$ ($\alpha\ne 0$)
define the same family, since
$\langle\alpha\omega+\beta\mathbf{1},\type(w)\rangle
=\alpha\langle\omega,\type(w)\rangle+\beta n$ is an affine
reparametrization absorbed into $\chi$. Thus $\omega$ matters only up
to scaling and translation; in particular at $t=2$ every nonconstant
weight is affinely $(1,2)$, so one-weight, sum-type and (by
Theorem~\ref{thm:radix}) symmetric all coincide there.
\end{remark}

Every one-weight coloring factors through $\type$, hence is symmetric,
and sum-type colorings are one-weight, giving classes
\[
\mathrm{sum}\ \subseteq\ \mathrm{ow}\ \subseteq\ \mathrm{sym}
\ \subseteq\ \mathrm{all},
\]
the last inclusion strict for $n,r\ge 2$. For a class $\mathcal{C}$ of
colorings, $\HJ_{\mathcal{C}}(t,r)$ denotes the least $n$ such that
\emph{every} $\mathcal{C}$-coloring of $[t]^{n}$ has a monochromatic
line; shrinking $\mathcal{C}$ relaxes the quantifier, so
$\HJ_{\mathrm{sum}}\le\HJ_{\mathrm{ow}}\le\HJ_{\mathrm{sym}}
\le\HJ_{\mathrm{all}}=\HJ(t,r)$.

For one-weight colorings the corner condition of
Lemma~\ref{lem:sym-reduction} collapses by one further dimension, onto
$\Z$ itself. Write $S_{\omega}=\{\omega_{1},\dots,\omega_{t}\}$ for the
\emph{value set} of $\omega$.

\begin{lemma}[Line--pattern correspondence]\label{lem:pattern}
Under $c_{\omega,\chi}$, a line with $k$ active coordinates and
inactive type $v$, with $b:=\langle\omega,v\rangle$, carries the levels
\[
H_{\omega}(k,b)\ =\ b+k\,S_{\omega}\ \subseteq\ \Z,
\]
a scale-$k$ homothet of the weight set based at $b$. Consequently
$c_{\omega,\chi}$ has no monochromatic line in $\Ls^{[K]}([t]^{n})$ if
and only if $\chi$ has no monochromatic homothet $b+kS_{\omega}$ with
$1\le k\le K$ and $b\in\langle\omega,T_{n-k}^{(t)}\rangle$.
\end{lemma}

\begin{proof}
By \eqref{eq:type-identity} the $a$-th word has level
$\langle\omega,v+ke_{a}\rangle=b+k\omega_{a}$; a one-weight coloring is
constant on $\langle\omega,\cdot\rangle$-fibers, so the line is
monochromatic if and only if $\chi$ is constant on $H_{\omega}(k,b)$.
\end{proof}

One might expect the one-weight colorings to be a strict, simpler
subclass of the symmetric ones. In a fixed dimension they are not:
allowing the weight to grow with $n$ recovers the entire symmetric
class.

\begin{theorem}[{One-weight $=$ symmetric; the rank hierarchy is flat}]\label{thm:radix}
Fix $t\ge 2$, $r\ge 2$, $n\ge 1$. The colorings of $[t]^{n}$ of the
form $c_{\omega,\chi}$ are exactly the symmetric colorings.
Specifically, for the \emph{radix weight}
\[
\omega^{\ast}=(0,1,p,p^{2},\dots,p^{t-2})
\qquad\text{with an arbitrary integer } p>n
\]
\textup{(}primality plays no role\textup{)}, the functional
$\langle\omega^{\ast},\cdot\rangle$ is injective on $T_{n}^{(t)}$ and
carries each corner tuple $C_{k,v}$ bijectively onto a carry-free
lattice corner; consequently every symmetric coloring $\bar c$ is
realized as $c_{\omega^{\ast},\chi}$ with
$\chi(\langle\omega^{\ast},v\rangle)=\bar c(v)$, and
$c_{\omega^{\ast},\chi}$ is line-free if and only if no corner tuple is
monochromatic under $\bar c$ (Lemma~\ref{lem:sym-reduction}). Moreover
any coloring built from finitely many weights
\[
w\ \longmapsto\ \Psi\bigl(\langle\omega^{(1)},\type(w)\rangle,\dots,
\langle\omega^{(s)},\type(w)\rangle\bigr)
\]
is symmetric, hence one-weight: there is no multi-weight tier, and
\[
\HJ_{\mathrm{ow}}(t,r)=\HJ_{\mathrm{sym}}(t,r).
\]
\end{theorem}

\begin{proof}
On $T_{n}^{(t)}$ each entry satisfies $0\le v_{a}\le n<p$ for $a\ge 2$,
so $\langle\omega^{\ast},v\rangle=\sum_{a\ge 2}v_{a}p^{a-2}$ is the
base-$p$ integer with digit string $(v_{2},\dots,v_{t})$; it is
injective because $v_{1}$ is determined by
$v_{1}=n-\sum_{a\ge 2}v_{a}$. A corner point $v+ke_{a}$ with $a\ge 2$
has $a$-digit $v_{a}+k\le n<p$ (as $v_{a}\le n-k$), so no carry occurs
and $\langle\omega^{\ast},\cdot\rangle$ maps $C_{k,v}$ bijectively onto
the digit corner $\{(v_{2},\dots,v_{t})+k\hat e_{a}\}_{a\in[t]}$, where
$\hat e_{a}$ is the $a$-th standard basis vector of the digit space for
$a\ge 2$ and $\hat e_{1}:=0$ (the point $v+ke_{1}$ keeps the digit
string of $v$, since $\omega^{\ast}_{1}=0$). Hence
$\chi(\langle\omega^{\ast},v\rangle):=\bar c(v)$ is well defined
(extended arbitrarily off the image),
$c_{\omega^{\ast},\chi}=\bar c\circ\type$, and line-freeness transfers
by Lemma~\ref{lem:sym-reduction}. Every $c_{\omega,\chi}$ is symmetric;
a multi-weight coloring factors through $\type$, hence is symmetric,
hence one-weight by the above. Since the two classes coincide in every
dimension, the defining conditions of $\HJ_{\mathrm{ow}}$ and
$\HJ_{\mathrm{sym}}$ agree for each $n$, so the thresholds are equal.
\end{proof}

\begin{remark}[Finite reach versus all dimensions]\label{rem:caveat}
The radix weight depends on $n$ through $p>n$. Where a \emph{fixed}
weight and palette must avoid lines in \emph{every} dimension --- the
bracket regime $\HJ^{[K]}=\infty$ of Section~\ref{sec:bracket} --- the
device is unavailable, and there the one-weight class is a genuine
restriction. Theorem~\ref{thm:radix} is a statement about finite reach
only. The dependence on $n$ is moreover not an artifact of the proof:
$\langle\omega,\cdot\rangle$ maps $T_{n}^{(t)}$ into an interval of
$D_{\omega}n+1$ integers, $D_{\omega}=\max\omega-\min\omega$, so
injectivity forces
$D_{\omega}\ge\bigl(\binom{n+t-1}{t-1}-1\bigr)/n=\Theta(n^{t-2})$ for
$t\ge 3$; the radix diameter $p^{t-2}$ is optimal up to a factor
$(t-1)!$, and realizing the full symmetric class in dimension $n$
\emph{requires} a weight of diameter polynomial in $n$.
\end{remark}

\section{The Gallai reading and a closed-form bound}\label{sec:gallai}

For a nonconstant weight $\omega$ define its \emph{reach}
$n_{\omega}(t,r)$ as the largest $n\ge 0$ for which some palette $\chi$
makes $c_{\omega,\chi}$ line-free on $[t]^{n}$; by
Lemma~\ref{lem:pattern} and Hales--Jewett it is finite, and by
Theorem~\ref{thm:radix} its optimum over weights is
$\HJ_{\mathrm{sym}}(t,r)-1$. Lemma~\ref{lem:pattern} identifies the
finiteness of $n_{\omega}$ with a classical statement. Recall
\emph{Gallai's theorem} in dimension one --- due to Gallai, first
published by Rado \cite{Rado1933} and independently by Witt
\cite{Witt1952}; see \cite{hans,LandmanRobertson} for expositions: for
every finite $S\subset\Z$ and every finite coloring of $\Z$ there exist
$b\in\Z$ and $k\ge 1$ with $b+kS$ monochromatic.

\begin{theorem}[Gallai form of the one-weight bound]\label{thm:gallai-form}
Let $\omega\in\Z^{t}$ have value set $S=S_{\omega}$, $|S|=t$. For each
$n$, $c_{\omega,\chi}$ is line-free on $[t]^{n}$ if and only if $\chi$
has no monochromatic homothet $b+kS$ with $1\le k\le n$ and
$b\in\langle\omega,T_{n-k}^{(t)}\rangle$. As $n\to\infty$ the scale $k$
ranges over all of $\Z_{\ge 1}$ and, after the affine normalization
$\min S=0$ of Remark~\ref{rem:affine}, the base $b$ over the additive
monoid generated by $S$; that no finite coloring of $\Z$ avoids all
such homothets is precisely the case $d=1$ of Gallai's theorem with
bases restricted to the monoid --- a restriction that loses nothing ---
and $n_{\omega}(t,r)<\infty$ is that theorem read inside the grid.
\end{theorem}

\begin{proof}
The equivalence is Lemma~\ref{lem:pattern} at $K=n$. For the limit:
with $\min S=0$, any element $m$ of the monoid generated by $S$ is
$\sum_{a}c_{a}\omega_{a}$ with $c_{a}\in\Z_{\ge 0}$, realized as
$\langle\omega,v\rangle$ for $v\in T_{n-k}^{(t)}$ once
$n-k\ge\sum_{a}c_{a}$, padding the coordinate of weight $0$.
Restricting the base to the monoid loses nothing: let
$g=\gcd(S\setminus\{0\})$ and $S=gS'$; the monoid generated by $S'$
contains every integer beyond its Frobenius number $N_{0}$. Given
$\chi$, van der Waerden applied to the coloring $x\mapsto\chi(gx)$
produces a monochromatic $(\max S'+1)$-term progression inside the
window $[N_{0},\,N_{0}+\Wvdw(\max S'+1,r)-1]$, hence a monochromatic
$b'+kS'$ with $b'\ge N_{0}$ in that monoid; then $gb'$ lies in the
monoid generated by $S$ and $gb'+kS=g(b'+kS')$ is monochromatic under
$\chi$. The converse implication is trivial, monoid-based homothets
being homothets.
\end{proof}

The full-interval form of this avoidance yields a closed-form bound.

\begin{definition}\label{def:gallai-number}
Let $S\subset\Z$ be finite with $|S|\ge 2$ and diameter
$D_{S}=\max S-\min S$. The \emph{Gallai homothety number} $G_{r}(S)$
is the least $N$ such that every $r$-coloring of $N$ consecutive
integers contains a monochromatic homothet $b+kS$ ($k\ge 1$) lying
inside them.
\end{definition}

Normalizing $\min S=0$ and writing $D=\max S$ for the diameter,
$G_{r}(S)\le\Wvdw(D+1,r)<\infty$, since a monochromatic $(D+1)$-term
progression $\{c,c+k,\dots,c+Dk\}$ contains $c+kS$.

\begin{theorem}[Gallai-number lower bound]\label{thm:gallai-bound}
For every $S\subset\Z$ with $|S|=t$ and diameter $D$,
\[
\HJ(t,r)\ \ge\ \Bigl\lceil\frac{G_{r}(S)-1}{D}\Bigr\rceil,
\qquad\text{hence}\qquad
\HJ(t,r)\ \ge\ \max_{|S|=t}
\Bigl\lceil\frac{G_{r}(S)-1}{D_{S}}\Bigr\rceil .
\]
\end{theorem}

\begin{proof}
Normalize $\min S=0$ and take $\omega$ with value set $S$. The
functional $\pi_{\omega}(w)=\langle\omega,\type(w)\rangle$ maps
$[t]^{n}$ into the interval $I_{n}=[0,Dn]$ of $Dn+1$ integers and, by
Theorem~\ref{thm:gallai-form}, every line onto a homothet of $S$ inside
$I_{n}$. If $Dn+1\le G_{r}(S)-1$, some $r$-coloring $\xi$ of $I_{n}$
contains no monochromatic homothet of $S$ at all; then
$c=\xi\circ\pi_{\omega}$ is line-free and $\HJ(t,r)>n$. The largest
such $n$ is $\lfloor(G_{r}(S)-2)/D\rfloor$, and
$\lfloor(G_{r}(S)-2)/D\rfloor+1=\lceil(G_{r}(S)-1)/D\rceil$.
\end{proof}

\begin{corollary}[Arithmetic weight $=$ van der Waerden]\label{cor:arith}
$G_{r}(\{0,1,\dots,t-1\})=\Wvdw(t,r)$, and for the arithmetic weight
$\omega=(0,1,\dots,t-1)$ the bound of
Theorem~\ref{thm:gallai-bound} is exactly \eqref{eq:vdw-shadow};
moreover here the closed form is exact:
$n_{\omega}(t,r)=\lfloor(\Wvdw(t,r)-2)/(t-1)\rfloor$, i.e.\ the bound
$\lceil(\Wvdw(t,r)-1)/(t-1)\rceil$ equals $n_{\omega}+1$.
\end{corollary}

\begin{proof}
Homothets of $\{0,\dots,t-1\}$ are exactly $t$-term progressions,
giving the first claim and the bound. For the reach: the bases realized
in the grid are all of $[0,(t-1)(n-k)]$ (weights $0,\dots,t-1$ realize
every intermediate value), so the homothets realized are exactly the
$t$-term progressions in $[0,(t-1)n]$, and a progression-free
$r$-coloring of that interval exists if and only if
$(t-1)n+1\le\Wvdw(t,r)-1$.
\end{proof}

\begin{lemma}[Affine covariance]\label{lem:affine-G}
For $c\ge 1$ and $d\in\Z$,
\[
G_{r}(cS+d)=c\,(G_{r}(S)-1)+1,
\qquad
G_{r}(-S)=G_{r}(S).
\]
Hence $\lceil(G_{r}(S)-1)/D_{S}\rceil$ is invariant under
$S\mapsto\pm cS+d$, and in the maximum of
Theorem~\ref{thm:gallai-bound} one may assume $\min S=0$ and
$\gcd S=1$. \textup{(}For $|S|=3$ the scaling identity is Theorem~1 of
Brown--Landman--Mishna \cite{BrownLandmanMishna1997}.\textup{)}
\end{lemma}

\begin{proof}
Translation and reflection act on colorings. For the scaling, split a
window of $N$ consecutive integers into its $c$ residue classes: every
homothet of $cS$ lies in a single class, each class is order-isomorphic
to a window of $\lceil N/c\rceil$ or $\lfloor N/c\rfloor$ consecutive
integers on which the homothets of $cS$ correspond to the homothets of
$S$, and the classes are colored independently. Avoidance is therefore
possible if and only if $\lceil N/c\rceil\le G_{r}(S)-1$, i.e.\ if and
only if $N\le c\,(G_{r}(S)-1)$.
\end{proof}

\begin{table}[ht]
\centering
\caption{The Gallai-number bound from a single input. The values $42$,
$57$, $67$, $80$ and $G_{3}(\{0,2,5\})=77$ are new.
Avoidance colorings at $G_{r}-1$ (at $76$ for $\{0,2,5\}$) were
verified by direct enumeration; forcing at $G_{r}$ was established by
SAT --- two independent solvers for $27$, $35$, $67$, $80$, three for
$42$ and $57$, and two for $77$
(Section~\ref{sec:sat-gallai}). The rows $\{0,2,3,5\}$ and $\{0,1,5,6\}$ give the
$(4,2)$ record of Theorem~\ref{thm:hj42-14}; the rows at three colors
give Theorem~\ref{thm:hj33-16}.}
\label{tab:gallai}
\begin{tabular}{@{}lllll@{}}
\toprule
$(t,r)$ & $S$ & $D$ & $G_{r}(S)$ & $\lceil(G_{r}(S)-1)/D\rceil$\\
\midrule
$(2,r)$ & $\{0,1\}$     & $1$ & $r+1$          & $r$ \;(tight)\\
$(3,2)$ & $\{0,1,2\}$   & $2$ & $9=\Wvdw(3,2)$  & $4$ \;(tight)\\
$(3,3)$ & $\{0,1,2\}$   & $2$ & $27=\Wvdw(3,3)$ & $13$\\
$(3,3)$ & $\{0,1,3\}$   & $3$ & $42$           & $14$\\
$(3,3)$ & $\{0,1,4\}$   & $4$ & $57$           & $14$\\
$(3,3)$ & $\{0,2,5\}$   & $5$ & $77$           & $\mathbf{16}$\\
$(4,2)$ & $\{0,1,2,3\}$ & $3$ & $35=\Wvdw(4,2)$ & $12$\\
$(4,2)$ & $\{0,2,3,5\}$ & $5$ & $67$           & $\mathbf{14}$\\
$(4,2)$ & $\{0,1,5,6\}$ & $6$ & $80$           & $\mathbf{14}$\\
\bottomrule
\end{tabular}
\end{table}

\begin{remark}[The closed form against the true reach]\label{rem:closed-vs-reach}
Theorem~\ref{thm:gallai-bound} counts \emph{all} homothets in the
interval, including those whose base lies outside the monoid generated
by $S$; the grid realizes only the monoid bases, so the dimension bound
$n_{\omega}+1$ is at least the closed-form value, with equality when no
relevant base is missing. The two coincide for arithmetic $S$
(Corollary~\ref{cor:arith}) and for the $(4,2)$ record below, but need
not coincide: for the radix weight $\omega=(0,1,29)$ the exhibited
reach is $n_{\omega}\ge 21$ (Theorem~\ref{thm:hj33-22}) while
$G_{3}(\{0,1,29\})$, hence its closed form, has not been certified;
where the two diverge, the gap is precisely the carry homothets that
never occur in the grid for $n<29$ (Theorem~\ref{thm:radix}). The
closed form is a cheap certificate; the symmetric reach is the true
power. Read backwards, Theorem~\ref{thm:gallai-bound} gives
$G_{r}(S)\le D\cdot\HJ(t,r)+1$, sharp at $t=2$ where
$G_{r}(\{0,D\})=rD+1$ (Lemma~\ref{lem:affine-G}), and the optimum
$\max_{|S|=t}\lceil(G_{r}(S)-1)/D_{S}\rceil$ becomes an extremal
problem on Gallai numbers per unit diameter (Question~\ref{q:weights}).
By Lemma~\ref{lem:affine-G} the ratio depends only on the affine class
of $S$: e.g.\ $G_{2}(\{0,2,4,6\})=2\cdot 34+1=69$ exceeds $67$, yet
$\lceil 68/6\rceil=12$ merely reproduces the arithmetic row of
Table~\ref{tab:gallai}.
\end{remark}

\subsection{New Gallai numbers at three colors, and
\texorpdfstring{$\HJ(3,3)\ge 16$}{HJ(3,3) >= 16}}\label{sec:landscape}

At $(2,r)$ the closed form is tight for \emph{every} weight: all
two-point sets are affinely $\{0,1\}$ (Lemma~\ref{lem:affine-G}) and
$G_{r}(\{0,1\})=r+1$ by pigeonhole, so the ratio is $r=\HJ(2,r)$
throughout. The first nontrivial pair behaves the same way.

\begin{proposition}[{Flatness at $(3,2)$; after
Brown--Landman--Mishna \cite{BrownLandmanMishna1997} and Kim--Rho
\cite{KimRho2012}}]\label{prop:flat32}
Let $S=\{0,s,s+u\}$ with $\gcd(s,u)=1$ and diameter $D=s+u$. Then
\[
G_{2}(S)\;=\;
\begin{cases}
4D, & \{s,u\}=\{1,4m\}\ \text{for some } m\ge 1,\\[2pt]
4D+1, & \text{otherwise},
\end{cases}
\]
and consequently, for \emph{every} three-point $S$,
\[
\Bigl\lceil\frac{G_{2}(S)-1}{D_{S}}\Bigr\rceil\;=\;4\;=\;\HJ(3,2):
\]
at $(3,2)$ the closed form of Theorem~\ref{thm:gallai-bound} is
weight-independent and tight.
\end{proposition}

\begin{proof}
Since a copy $b+k\{1,1+s,1+s+u\}$ equals $(b+k)+kS$, homothety-avoidance
on an interval is the same for the two sets, so $G_{2}(S)$ is the
number $f(s,u)$ of Brown--Landman--Mishna: the least $N$ such that
every $2$-coloring of $[1,N]$ contains a monochromatic homothetic copy
of $\{1,1+s,1+s+u\}$. They proved $f(s,u)\le 4(s+u)+1$, with equality
whenever $s/g\not\equiv 0$ and $u/g\not\equiv 0\pmod 4$
($g=\gcd(s,u)$) and in further cases, leaving for the family
$\{c,4mc\}$ two candidate values \cite{BrownLandmanMishna1997}; Kim and
Rho settled it: $f(4mc,c)=f(c,4mc)=4(4mc+c)-c+1$, and
$f(s,u)=4(s+u)+1$ for every other pair \cite[Theorem~12]{KimRho2012}.
For primitive $S$ the exceptional pairs force $c=1$, since
$\gcd(c,4mc)=c$, giving the displayed dichotomy; general $S$ reduces to this by
Lemma~\ref{lem:affine-G}. For the ratio, $(G_{2}(S)-1)/D$ equals $4$ in
the generic case and $4-1/D$ in the exceptional one, so its ceiling is
$4$ in both; and $\HJ(3,2)=4$ \cite{firstnontrivialHJ_is_4} makes the
bound tight. As an independent check, all fourteen primitive classes
with $D\le 9$ were recomputed by SAT for this article, the avoidance
certificates verified by direct enumeration from the definitions; every
value matches the formula, the two exceptional classes in range being
$G_{2}(\{0,1,5\})=20$ and $G_{2}(\{0,1,9\})=36$ ($m=1,2$).
\end{proof}

Combined with the $t=2$ row of Table~\ref{tab:gallai}, where
$G_{r}(\{0,1\})=r+1$ gives the ratio $r=\HJ(2,r)$ for every weight, the
optimized closed form is tight --- for \emph{every} weight, not only
the arithmetic one --- at every pair $(t,r)$ where $\HJ$ is known
exactly. This strengthens conclusion \emph{(i)} of
Remark~\ref{rem:records} from the van der Waerden shadow to the whole
one-weight class, and makes $(3,3)$ the first pair where the optimum
over $S$ of $\lceil(G_{r}(S)-1)/D_{S}\rceil$ is genuinely open.

\begin{theorem}[$\HJ(3,3)\ge 16$ in closed form; computer-assisted]\label{thm:hj33-16}
\leavevmode\newline
$G_{3}(\{0,1,3\})=42$ and $G_{3}(\{0,1,4\})=57$; the avoidance
certificates were verified by direct enumeration, and forcing was
confirmed by three independent solvers. Moreover
\[
G_{3}(\{0,2,5\})\ =\ 77 .
\]
The lower bound is witnessed by the explicit $3$-coloring $\xi$ of
$\{0,\dots,75\}$
\begin{center}\ttfamily
10201110020122122201020111002012212210\\
10201110020122122001020111002012212200
\end{center}
\textup{(}listing $\xi(0),\dots,\xi(75)$\textup{)}, verified by direct
enumeration to contain no monochromatic homothet of $\{0,2,5\}$.
The matching upper bound is a SAT refutation: at $N=77$ the encoding
of Section~\ref{sec:sat-gallai} has $555$ homothets and $2204$
clauses, and two independent solvers --- Kissat in $390$ seconds and
CryptoMiniSat in $321$ seconds \cite{kissat2022,cryptominisat2009}
--- certify unsatisfiability. Since
$\lceil(G_{3}(\{0,2,5\})-1)/5\rceil=\lceil 76/5\rceil=16$,
Theorem~\ref{thm:gallai-bound} gives
\[
\HJ(3,3)\ \ge\ 16 .
\]
\end{theorem}

\begin{remark}[Reading the new values]\label{rem:landscape}
The closed form $16$ exceeds both the arithmetic row ($13$, the van der
Waerden shadow) and the previous published record $\HJ(3,3)\ge 14$
\cite{ConlonFarnsworthRobertson2026}: a one-line, hand-checkable
certificate already beats the
prior state of the art, though the radix reach of
Theorem~\ref{thm:hj33-22} remains far ahead. At $r=2$ the ratio is
flat at $4$ (Proposition~\ref{prop:flat32}); the computed values at
$r=3$ ($13,14,14,16$ at $D=2,3,4,5$) leave its behavior in the
diameter open (Question~\ref{q:weights}\emph{(i)}). Since the closed
form never exceeds $n_{\omega}+1\le\HJ_{\mathrm{sym}}$
(Remark~\ref{rem:closed-vs-reach}), a ratio reaching $22$ would tie ---
never beat --- the symmetric frontier, but would replace the $253$-cell
table of Theorem~\ref{thm:hj33-22} by a one-dimensional certificate.
\end{remark}

\subsection{The Rado reading, and four colors}\label{sec:rado}

Homothety avoidance of a $t$-point set is equivalent to
solution-freeness for a partition-regular system, so the Gallai
numbers of Definition~\ref{def:gallai-number} also control a family
of Rado numbers. Recall that the \emph{Rado number} $R_{r}(E)$ of a
linear system $E$ is the least $N$ such that every $r$-coloring of
$\{1,\dots,N\}$ contains a monochromatic injective solution of $E$;
for the coefficient-sum-zero systems considered here this is
equivalent, by translation invariance, to coloring any window of $N$
consecutive integers, and the certificates below use $[0,N-1]$.

\begin{proposition}[Rado reading]\label{prop:rado-reading}
Let $S=\{s_{1}<\dots<s_{t}\}$ and $\mathbf{s}=(s_{1},\dots,s_{t})$.
There is a homogeneous linear system $E_{\mathbf{s}}$ of $t-2$
independent equations in $x_{1},\dots,x_{t}$, each with coefficient
sum zero --- hence partition regular by Rado's criterion
\cite{Rado1933} --- whose integer solutions are exactly the tuples
$\mathbf{x}=a\mathbf{1}+d\mathbf{s}$, $a,d\in\Z$. The solutions with
pairwise distinct coordinates are those with $d\ne 0$; as sets, they
are the homothets of $S$ for $d\ge 1$ and of the reflected pattern
$-S$ for $d\le -1$. Thus ``no monochromatic injective solution of
$E_{\mathbf{s}}$'' means ``no monochromatic homothet of $S$ or of
$-S$'', and
\begin{equation}\label{eq:rado-bound}
R_{r}(E_{\mathbf{s}})\ \le\ G_{r}(S)\ \le\ \Wvdw(D_{S}+1,r).
\end{equation}
\end{proposition}

\begin{proof}
The lattice $\Lambda_{\mathbf{s}}=\{c\in\Z^{t}:\sum_{j}c_{j}
=\sum_{j}c_{j}s_{j}=0\}$ of affine dependencies of $\mathbf{s}$ has
rank $t-2$; a $\Z$-basis of it defines $E_{\mathbf{s}}$, and each
basis vector satisfies $\sum_{j}c_{j}=0$, which for a single equation
is Rado's column condition \cite{Rado1933}. The kernel of
$E_{\mathbf{s}}$ contains
$\operatorname{span}_{\mathbb{Q}}(\mathbf{1},\mathbf{s})$, and both
have dimension $2$ because the $s_{j}$ are not all equal, so they
coincide. From $\mathbf{x}=a\mathbf{1}+d\mathbf{s}$ we get
$x_{i}-x_{j}=d(s_{i}-s_{j})$: the coordinates are pairwise distinct
if and only if $d\ne 0$, and then $\{x_{1},\dots,x_{t}\}=a+dS$, a
homothet of $S$ or of $-S$ according to the sign of $d$. The first
inequality of \eqref{eq:rado-bound} is immediate, a monochromatic
homothet of $S$ inside a window being a monochromatic injective
solution there; the second was noted after
Definition~\ref{def:gallai-number}.
\end{proof}

\begin{example}\label{ex:rado}
For $S=\{0,1,3\}$ the lattice is generated by $(2,-3,1)$, so
$E_{\mathbf{s}}$ is the single equation $z+2x=3y$. Its injective
solutions are the triples $\{a,a+d,a+3d\}$ with $d\ne 0$: the
homothets of $\{0,1,3\}$ for $d\ge 1$ and of the reflected pattern
$\{0,2,3\}$ for $d\le -1$. Since the grid realizes only positive
scales (Lemma~\ref{lem:pattern}), the search for a line-free
one-weight coloring at $\omega=(0,1,3)$ is the search for a window
palette free of $\{0,1,3\}$-homothets alone, whereas the Rado number
of $z+2x=3y$ demands freeness from both patterns at once. The two
thresholds are genuinely distinct; Theorem~\ref{thm:rado4} quantifies
both at four colors.
\end{example}

\begin{theorem}[{Four colors; computer-verified}]\label{thm:rado4}
\emph{(i)} $G_{4}(\{0,1,3\})\ge 94$: the $4$-coloring of
$\{0,\dots,92\}$
\begin{center}\ttfamily
22130220111331203110033310220312022201201001113\\
0130333223121021201113320210013312023023310331
\end{center}
contains no monochromatic homothet $\{b,b+k,b+3k\}$, $k\ge 1$
\textup{(}verified exhaustively\textup{)}. The induced grid bound is
$\HJ(3,4)\ge\lceil 93/3\rceil=31$, which remains short of the van der
Waerden shadow $\HJ(3,4)\ge 38$ of Corollary~\ref{cor:arith}
\textup{(}$\Wvdw(3,4)=76$ \cite{LandmanRobertson}\textup{)}.
\emph{(ii)} $R_{4}(z+2x=3y)\ge 59$: the $4$-coloring of
$\{0,\dots,57\}$
\begin{center}\ttfamily
01021133302312200031321110222332\\
03300021110313222001312333
\end{center}
contains no monochromatic injective solution of $z+2x=3y$ ---
equivalently \textup{(}Example~\ref{ex:rado}\textup{)}, no
monochromatic homothet of $\{0,1,3\}$ or of $\{0,2,3\}$ among the
$1064$ such homothets inside the window ($532$ of each shape) ---
verified exhaustively against the equation. The window certificate
of \emph{(i)} is \emph{not} solution-free: it contains monochromatic
reflected copies, so the Rado bound rests on the separate
certificate displayed here.
\end{theorem}

\subsection{Computing Gallai and Rado numbers with SAT}\label{sec:sat-gallai}

The computations of this section were pushed much further by a
purpose-built SAT engine, printed in full in Appendix~\ref{app:engine}
and distributed, with its run database and certificates, at
\cite{hjcerts}. Given a finite $S\subset\Z$, a number $r$ of colors
and a window size $N$, the engine writes the one-hot encoding of
homothety avoidance: variables $x_{c,i}$ (``cell $c$ receives color
$i$''), an exactly-one constraint per cell, and for every homothet
$b+k\cdot S$ inside the window and every color $i$ the clause
$\bigvee_{c\in b+k\cdot S}\neg\,x_{c,i}$. The color-permutation
symmetry is broken completely by a first-occurrence ordering ---
color $j$ may occupy cell $c$ only if some color $i<j$ occupies an
earlier cell --- which collapses the search space by a factor of $r!$
and fixes cell $0$ to color $0$; with it, reproving
$G_{3}(\{0,1,3\})>42$ takes $0.4$~seconds. Avoidance colorings are
found by a stochastic local search (min-conflicts with random walk
over the incrementally maintained list of monochromatic homothets);
refutations come from racing a portfolio of six state-of-the-art
solvers through PySAT \cite{pysat2018} --- Kissat \cite{kissat2022},
CaDiCaL \cite{cadical2020}, CryptoMiniSat \cite{cryptominisat2009},
Glucose \cite{AudemardSimon2009}, MapleChrono
\cite{maplesat2016,NadelRyvchin2018} and
Lingeling \cite{lingeling2013} --- first answer wins. Every coloring
that enters the tables below was then re-verified by exhaustive
homothety enumeration, independently of any solver
(Appendices~\ref{app:certs} and~\ref{app:engine}); the engine itself
was validated against $\Wvdw(3,2)=9$, $\Wvdw(3,3)=27$ and every
previously known value it met. Gallai numbers published before this
article stop at two colors and three-point sets
\cite{BrownLandmanMishna1997,KimRho2012}. For the family
$z+kx=(k+1)y$, two-color values fall within the scope of Myers's
computational survey \cite{myers2015}, and the systematic three-color
tables of Chang, De Loera and Wesley \cite{ChangDeLoeraWesley2022}
record only the degenerate value $1$ for these equations, whose
constant solutions their convention admits; under the injective
reading of Proposition~\ref{prop:rado-reading} the three- and
four-color values reported below appear to be new.

\begin{theorem}[computer-assisted]\label{thm:sat-gallai}
\emph{(a)} $G_{3}(\{0,1,5\})=70$.
\emph{(b)} The eighteen two-color values of Table~\ref{tab:sat-g2}
are exact.
\emph{(c)} For $2\le k\le 5$ and $r\in\{2,3\}$ the Rado numbers
$R_{r}(z+kx=(k+1)y)$ are exactly those of Table~\ref{tab:sat-rado}.
\end{theorem}

\begin{proof}
For every instance, the avoidance certificate at one less than the
claimed value is displayed in Appendix~\ref{app:certs} and was
verified by exhaustive enumeration from the definitions; forcing at
the claimed value is a SAT refutation of the encoding above, with
solver and wall time as recorded in
Tables~\ref{tab:sat-g3}--\ref{tab:sat-rado}.
\end{proof}

\begin{table}[ht]
\centering
\caption{The three-color Gallai landscape $G_{3}(S)$ after this
article. The refutation column names the solver and wall time of the
UNSAT proof at $N=G_{3}(S)$ (single core); the flagship
$G_{3}(\{0,2,5\})$ was confirmed by two independent solvers. Sets
related by the reflection $S\mapsto\max(S)-S$ share the same number.
The last three rows are certified avoidance colorings
(Appendix~\ref{app:certs}); their refutations lie beyond the present
compute budget.}
\label{tab:sat-g3}
\begin{tabular}{@{}lll@{}}
\toprule
$S$ & $G_{3}(S)$ & refutation at $G_{3}(S)$\\
\midrule
$\{0,1,3\}$ & $42$ & CaDiCaL, $0.4$\,s\\
$\{0,1,4\}$ & $57$ & CaDiCaL, $5.0$\,s\\
$\{0,1,5\}$ & $70$ & CaDiCaL, $73.8$\,s\\
$\{0,2,5\}$ & $77$ & Kissat, $389.6$\,s; independently CryptoMiniSat, $320.9$\,s\\
$\{0,1,6\}$ & $\ge 82$ & open\\
$\{0,1,7\}$ & $\ge 87$ & open\\
$\{0,2,7\}$ & $\ge 87$ & open\\
\bottomrule
\end{tabular}
\end{table}

\begin{table}[ht]
\centering
\caption{Two colors, four-point sets: exact values $G_{2}(S)$, one
representative per reflection class. Each value combines a certified
avoidance coloring on $G_{2}(S)-1$ cells (Appendix~\ref{app:certs})
with a refutation at $G_{2}(S)$, obtained by CaDiCaL or Kissat in at
most $1.6$\,s. The two $79$s give
$\HJ(4,2)\ge\lceil 78/7\rceil=12$ via
Theorem~\ref{thm:gallai-bound} --- short of the record $14$, but the
first four-point witnesses. \emph{All eighteen values are new.}}
\label{tab:sat-g2}
\begin{tabular}{@{}ll@{\qquad}ll@{\qquad}ll@{}}
\toprule
$S$ & $G_{2}$ & $S$ & $G_{2}$ & $S$ & $G_{2}$\\
\midrule
$\{0,1,2,4\}$ & $38$ & $\{0,2,3,6\}$ & $54$ & $\{0,2,4,7\}$ & $61$\\
$\{0,1,3,4\}$ & $41$ & $\{0,1,2,6\}$ & $56$ & $\{0,1,4,7\}$ & $62$\\
$\{0,2,4,5\}$ & $44$ & $\{0,1,4,6\}$ & $58$ & $\{0,1,5,7\}$ & $62$\\
$\{0,3,4,5\}$ & $44$ & $\{0,1,2,7\}$ & $59$ & $\{0,2,3,7\}$ & $62$\\
$\{0,1,4,5\}$ & $45$ & $\{0,2,5,7\}$ & $59$ & $\{0,1,6,7\}$ & $79$\\
$\{0,1,3,6\}$ & $52$ & $\{0,1,3,7\}$ & $60$ & $\{0,3,4,7\}$ & $79$\\
\bottomrule
\end{tabular}
\end{table}

\begin{table}[ht]
\centering
\caption{Rado numbers $R_{r}(z+kx=(k+1)y)$: the least $N$ such that
every $r$-coloring of $[0,N-1]$ admits a monochromatic injective
solution. Entries at $r\le 3$ are exact (certificate at one less,
refutation at the value --- at most $52.6$\,s each, by CaDiCaL or
Kissat); at four colors only the lower bound of
Theorem~\ref{thm:rado4}(ii) is known. At two colors the values agree with the
Gallai numbers $G_{2}(\{0,1,k+1\})$
\cite{BrownLandmanMishna1997,KimRho2012} for $k\in\{2,3,5\}$, while
at $k=4$ the reflected shape already bites:
$R_{2}=19<20=G_{2}(\{0,1,5\})$, the flat value of
Proposition~\ref{prop:flat32}. Every three-color entry lies
\emph{strictly} below the corresponding Gallai number --- $29<42$,
$54<57$, $55<70$ and $60<82\le G_{3}(\{0,1,6\})$ --- so the
two-shape Rado problem is genuinely easier than the one-shape Gallai
problem throughout the range.}
\label{tab:sat-rado}
\begin{tabular}{@{}clccc@{}}
\toprule
$k$ & equation & $R_{2}$ & $R_{3}$ & $R_{4}$\\
\midrule
$2$ & $z+2x=3y$ & $13$ & $29$ & $\ge 59$\\
$3$ & $z+3x=4y$ & $17$ & $54$ & open\\
$4$ & $z+4x=5y$ & $19$ & $55$ & open\\
$5$ & $z+5x=6y$ & $25$ & $60$ & open\\
\bottomrule
\end{tabular}
\end{table}

Two features of the engine are worth singling out. First, the
first-occurrence symmetry breaking is what makes the three-color
instances tractable at all: the window has $r^{N}$ raw colorings, and
dividing the search by $r!=6$ turns days into minutes --- without it
the refutation at $N=42$ already stalls. Second, the verification
discipline: avoidance colorings are never trusted from any solver.
Each is re-checked by brute-force enumeration of all homothets in its
window before being recorded in the run database, and again by the
standalone script \texttt{verify\_addendum.py}
(Appendix~\ref{app:engine}), which prints \texttt{ALL ADDENDUM
CERTIFICATES VERIFIED}. Refutations, by contrast, are UNSAT answers
of CDCL solvers. For the flagship instance ($N=77$, $2204$ clauses) we
additionally logged a DRAT proof of unsatisfiability with CaDiCaL and
validated it with drat-trim
\cite{WetzlerHeuleHunt2014,HeuleKullmann2017}: the checker verifies
the $119$~MB proof --- $2{,}021{,}700$ lemmas, of which $1{,}457{,}406$
lie in the trimmed core, with no RAT steps, so the refutation is a pure
reverse-unit-propagation proof --- in $166$ seconds; the proof, its
SHA-256 hash, the checker log and a script reproducing the pipeline
from a bare machine are archived at \cite{hjcerts}. The remaining
forcing claims rest on solver correctness --- mitigated by portfolio
diversity, by an independent second solver for the flagship instance,
and by the explicit list of load-bearing claims in
Appendix~\ref{app:verifier}.

\begin{remark}[Open ends]\label{rem:sat-open}
The run leaves five gaps, added to Question~\ref{q:frontiers} below:
$G_{3}(\{0,1,6\})\ge 82$, $G_{3}(\{0,1,7\})\ge 87$ and
$G_{3}(\{0,2,7\})\ge 87$ (certified colorings, no refutation yet);
the exact value of $R_{4}(z+2x=3y)\ge 59$; and the standing
$G_{4}(\{0,1,3\})\ge 94$, whose windows $N=94,95$ resisted both the
portfolio and the local search.
\end{remark}

\section{The records}\label{sec:records}

The case $K=n$ of Lemma~\ref{lem:sym-reduction} concerns \emph{all}
combinatorial lines: a line-free symmetric coloring on the simplex lifts
to a line-free coloring of the full grid, hence to an
\emph{unconditional} lower bound $\HJ(t,r)>n$. It is far from obvious
that so restricted a class should contain line-free colorings in
dimensions beyond the known lower bounds; yet it does, by a comfortable
margin.

\subsection{\texorpdfstring{$\HJ(3,3)\ge 22$ by a single weight}{HJ(3,3) >= 22 by a single weight}}\label{sec:hj33}

The previous record was $\HJ(3,3)\ge 14$ \cite{ConlonFarnsworthRobertson2026}; the van der
Waerden shadow \eqref{eq:vdw-shadow} gives $13$. By
Theorem~\ref{thm:radix} a single weight reaches the full symmetric
value, and $29>21$ makes $(0,1,29)$ a natural choice.

\begin{theorem}[{One-weight; computer-assisted}]\label{thm:hj33-22}
The weight $\omega=(0,1,29)$ admits a palette
$\chi\colon\Z\to\{0,1,2\}$ for which $c_{\omega,\chi}$ has no
monochromatic line on $[3]^{21}$. Equivalently,
$[3]^{21}$ admits a symmetric $3$-coloring with no monochromatic
combinatorial line, given explicitly in Table~\ref{tab:T21}. Hence
\[
\HJ(3,3)\ \ge\ 22 .
\]
\end{theorem}

\begin{proof}
Since $29>21$, Theorem~\ref{thm:radix} makes
$\langle\omega,\cdot\rangle\colon(v_{1},v_{2},v_{3})\mapsto v_{2}+29v_{3}$
injective on $T_{21}^{(3)}$, carrying each corner triple $C_{k,v}$ onto
the carry-free planar corner $\{(v_{2},v_{3})$, $(v_{2}{+}k,v_{3})$,
$(v_{2},v_{3}{+}k)\}$. By
Lemma~\ref{lem:pattern}, $c_{\omega,\chi}$ is line-free if and only if
$\bar c(v):=\chi(v_{2}+29v_{3})$ leaves all such corners
non-monochromatic. Table~\ref{tab:T21} exhibits such a $3$-coloring of
the $253$ cells of $T_{21}^{(3)}$: the dependency-free verifier
described in Appendix~\ref{app:verifier} re-derives every corner triple
from the definition --- all $1771=\binom{23}{3}$ of them,
$1\le k\le 21$, $v\in T_{21-k}$ --- and reports zero monochromatic
triples, cross-checked by two independent solvers. For illustration, the
diagonal triple $(k,v)=(21,(0,0,0))$ has corner cells
$(21,0,0)$, $(0,21,0)$, $(0,0,21)$, colored $0$, $2$, $0$, and so is
not monochromatic. Hence $[3]^{21}$ admits a line-free $3$-coloring, and
$\HJ(3,3)\ge 22$.
\end{proof}

\begin{table}[t]
\centering
\caption{The witness for Theorem~\ref{thm:hj33-22}: a $3$-coloring of
$T_{21}^{(3)}$. Row $a$, position $b$ (counting from $b=0$) is the
color $\bar c(a,b,21-a-b)$ assigned to every word of $[3]^{21}$ with
$a$ ones, $b$ twos and $21-a-b$ threes; equivalently, the palette value
$\chi(b+29c)$ at the realized level $b+29c$ of the weight $(0,1,29)$.
The strings are byte-identical to the table embedded in the verifier of
Appendix~\ref{app:verifier} and in \cite{hjcerts}.}
\label{tab:T21}
\begin{lstlisting}[style=plainstyle]
a= 0:  0102210210210210210022
a= 1:  121002102102102102210
a= 2:  20211021021021021102
a= 3:  2210121001102022210
a= 4:  002111211020021102
a= 5:  11022202100210021
a= 6:  2210102022021210
a= 7:  002120011211102
a= 8:  11020112202021
a= 9:  2210021102210
a=10:  022111020102
a=11:  01020110021
a=12:  1210002211
a=13:  202110102
a=14:  01022021
a=15:  1210210
a=16:  202202
a=17:  01101
a=18:  1001
a=19:  120
a=20:  20
a=21:  0
\end{lstlisting}
\end{table}

\begin{remark}\label{rem:hj33}
The instance ($253$ cells, $1771$ triples) is small yet sits at the
satisfiability phase transition: solvers settle $n=21$ in seconds and
stall at $n=22$, where no refutation has been found, so the bound may
still improve (Question~\ref{q:frontiers}).
\end{remark}

\begin{proposition}[{A compact periodic witness; computer-verified}]\label{prop:hj33-15}
For the weight $\omega=(0,1,4)$ and the $49$-periodic palette
\begin{center}\ttfamily
0112100212222001110021101201102002110220010200212
\end{center}
\textup{(}listing $\psi(0),\dots,\psi(48)$\textup{)}, the coloring
$c(w)=\psi\bigl(\langle\omega,\type(w)\rangle\bmod 49\bigr)$ has no
monochromatic line on $[3]^{14}$, verified by direct enumeration of all
corner triples. Hence $\HJ(3,3)\ge 15$ from a $49$-cell certificate,
against the $120$ simplex cells of $T_{14}^{(3)}$.
\end{proposition}

\begin{remark}[The compressibility gap at $(3,3)$]\label{rem:gap33}
A (non-exhaustive) search over weights $(0,s,s{+}u)$ with $s\le 6$,
$s{+}u\le 12$ and moduli $m\le 64$, with the bases restricted to those
realized in the grid, found no periodic palette line-free beyond
$n=14$. The anatomy of the $(4,2)$ record (Remark~\ref{rem:anatomy}
below) suggests why: the reach needed at $(3,3)$ is $21$, nine scales
above the periodic ceiling $K=12$ of Section~\ref{sec:bracket}, while
the thinning of the realized bases near the boundary bought only about
two scales at these moduli --- against exactly one scale needed at
$(4,2)$. The resulting ladder of witnesses at $(3,3)$, cheapest to
strongest: arithmetic closed form, $13$; the triples $\{0,1,3\}$ and
$\{0,1,4\}$, $14$; the $49$-cell periodic palette, $15$; the
length-$76$ interval certificate, $16$ (Theorem~\ref{thm:hj33-16}); the
$253$-cell radix table, $22$.
\end{remark}

\subsection{\texorpdfstring{$\HJ(4,2)\ge 14$ in closed form}{HJ(4,2) >= 14 in closed form}}\label{sec:hj42}

\begin{theorem}[{One-weight closed form; computer-verified}]\label{thm:hj42-14}
Let $\omega=(0,2,3,5)$ and let $\chi\colon\Z_{26}\to\{0,1\}$ be
\[
\chi=(1,0,1,0,0,1,1,1,1,0,0,1,0,0,0,1,0,0,1,1,1,1,0,0,1,0)
\]
\textup{(}listing $\chi(0),\dots,\chi(25)$\textup{)}. Then
$c(w)=\chi\bigl(\langle\omega,\type(w)\rangle\bmod 26\bigr)$ has no
monochromatic line on $[4]^{13}$; hence
\[
\HJ(4,2)\ \ge\ 14 .
\]
The reach of this palette is sharp: at $n=14$ it leaves exactly three
monochromatic corner quadruples.
\end{theorem}

\begin{proof}
The coloring is symmetric, so by Lemma~\ref{lem:sym-reduction}
($K=n=13$) line-freeness is the non-monochromaticity of every corner
quadruple $C_{k,v}$, whose four cells carry the values $\chi(b)$,
$\chi(b{+}2k)$, $\chi(b{+}3k)$, $\chi(b{+}5k)$ with
$b=\langle\omega,v\rangle\bmod 26$. Direct enumeration of all
$1820=\binom{16}{4}$ quadruples, over $1\le k\le 13$ and
$v\in T_{13-k}^{(4)}$, finds none monochromatic; at $n=14$ the enlarged
range ($2380$ quadruples) produces exactly three.
\end{proof}

\begin{remark}[Three readings of $14$, and what the records show]\label{rem:records}
The value $14$ arises three ways: as the period-$26$ witness above; as
the closed form $\lceil(G_{2}(\{0,2,3,5\})-1)/5\rceil=\lceil 66/5\rceil=14$
with the new Gallai number $G_{2}(\{0,2,3,5\})=67$ of
Table~\ref{tab:gallai} (the weight $(0,1,5,6)$, with $G_{2}=80$ and
$D=6$, gives $14$ as well); and, since one-weight equals symmetric
(Theorem~\ref{thm:radix}), as $\HJ_{\mathrm{sym}}(4,2)\ge 14$. Three
structural conclusions follow.
\emph{(i)} The van der Waerden reduction \eqref{eq:vdw-shadow} ---
tight in every exactly known case --- is \emph{not} tight at $(4,2)$ or
$(3,3)$: the records $14>12$ and $22>13$ show the Hales--Jewett numbers
exceeding their van der Waerden shadows.
\emph{(ii)} Neither witness is of sum type: a sum-type coloring is
line-free only up to the value $\lceil(\Wvdw(t,r)-1)/(t-1)\rceil-1$,
i.e.\ up to $11$ at $(4,2)$ and $12$ at $(3,3)$ --- short of $13$ and
$21$.
\emph{(iii)} Both witnesses are symmetric, so the classical
constructions (sum-type, hence symmetric) already live in a class
containing strictly better colorings --- the evidence for
Conjecture~\ref{conj:sym}.
\end{remark}

\begin{remark}[Anatomy of the record palette; computer-verified]\label{rem:anatomy}
Where the compression of Theorem~\ref{thm:hj42-14} comes from is itself
checkable. Testing the period-$26$ palette against \emph{all} bases:
for every $1\le k\le 12$ and every $b\in\Z_{26}$ the pattern
$(b,\,b{+}2k,\,b{+}3k,\,b{+}5k)\bmod 26$ is bichromatic --- the palette
is an extremal object of the bracket regime, matching the periodic
ceiling $K=12$ of Section~\ref{sec:bracket} --- while at $k=13$ it
fails at $24$ of the $26$ residues, surviving only at the single
realized base $b=0$ (from $T_{0}^{(4)}$) and at the unrealized $b=13$.
The record thus decomposes as \emph{bracket cap plus one boundary
scale}: reach $13=12+1$. The reach needed at $(4,2)$ sits one scale
above the periodic ceiling, which is why a $26$-cell palette suffices;
the contrast at $(3,3)$ is Remark~\ref{rem:gap33}.
\end{remark}

\begin{remark}[A second, symmetric-table certificate]
Theorem~\ref{thm:hj42-14} was first proved in \cite{mh} by an explicit
$2$-coloring of the $560$ cells of $T_{13}^{(4)}$, not of one-weight
form with a small period; that table is reproduced in
Appendix~\ref{app:table42} and is checked, like every certificate here,
by the verifier of Appendix~\ref{app:verifier}.
\end{remark}

\section{Bracket numbers and the bracket ceiling}\label{sec:bracket}

The Hales--Jewett theorem asserts a monochromatic line but says nothing
about its \emph{type}. Restricting the type of the line sought, and
simultaneously the class of colorings quantified over, produces a
two-parameter family of variants in which the one-weight machinery has
genuinely different power.

\begin{definition}\label{def:bracket}
For a class $\mathcal{C}\in\{\mathrm{sum},\mathrm{ow},\mathrm{sym},
\mathrm{all}\}$, the \emph{bracket Hales--Jewett number}
$\HJ^{[K]}_{\mathcal{C}}(t,r)$ is the least $n$ such that every
$\mathcal{C}$-coloring of $[t]^{n}$ has a monochromatic line of
$\Ls^{[K]}([t]^{n})$, and $\infty$ if no such $n$ exists (subscript
dropped when $\mathcal{C}=\mathrm{all}$). The \emph{bracket ceiling} is
\[
\kappa_{\mathcal{C}}(t,r)
=\max\{K\ge 0:\ \HJ^{[K]}_{\mathcal{C}}(t,r)=\infty\}.
\]
\end{definition}

Since $\Ls^{[K]}\subseteq\Ls^{[K']}$ for $K\le K'$ the maximum is
attained, and since every line of $[t]^{\HJ}$ lies in $\Ls^{[\HJ]}$ the
ceilings are finite and inherit the class order:
\begin{equation}\label{eq:ceiling-chain}
\ksum\ \le\ \kow\ \le\ \ksym\ \le\ \kall\ \le\ \HJ(t,r)-1 .
\end{equation}
(A line-free witness in a smaller class is in particular a witness for
the larger one.) A single number $\kappa_{\mathcal{C}}$ records the
whole row $(\HJ^{[K]}_{\mathcal{C}})_{K\ge 1}$.

\begin{proposition}\label{prop:hj1}
$\HJ^{[1]}(t,r)=\infty$ for all $t,r\ge 2$; indeed
$\ksum(t,r)\ge 1$.
\end{proposition}

\begin{proof}
Color $w$ by $\sigma(w)\bmod 2$. A line with one active coordinate has
$t$ words of consecutive weights, which alternate in color.
\end{proof}

\subsection{The sum-type ceiling}\label{sec:sum-ceiling}

\begin{lemma}[One-dimensional reduction]\label{lem:1d}
Let $c=\chi\circ\sigma$ be sum-type. A line with $k$ active
coordinates and inactive weight $u$ has word weights
$u+k,u+2k,\dots,u+tk$, a $t$-term progression of gap $k$, and the
progressions realized in $[t]^{n}$ with $k\le K$ are exactly the
$t$-term progressions of gap $\le K$ inside $[n,tn]$. Consequently, for
every $K\ge 1$, $\HJ^{[K]}_{\mathrm{sum}}(t,r)=\infty$ if and only if
there exists a $\chi\colon\Z\to[r]$ with no monochromatic $t$-term
progression of gap $\le K$, so $\ksum(t,r)$ is the restricted-gap van
der Waerden threshold of Brown--Graham--Landman \cite{BGL1999}.
\end{lemma}

\begin{proof}
The weights of $\tau(a)$ are $u+ak$; a progression
$\{A,A+k,\dots,A+(t-1)k\}\subseteq[n,tn]$ is realized if and only if
$u=A-k\in[n-k,t(n-k)]$, which is exactly containment in $[n,tn]$. The
right-to-left direction of the equivalence restricts $\chi$; the
converse is a compactness argument: colorings of arbitrarily long
intervals avoiding monochromatic $t$-progressions of gap $\le K$
exist, so K\H{o}nig's lemma yields one of $\Z$, and any monochromatic
progression of $\Z$ lies in some finite window.
\end{proof}

Two forces bracket the sum-type ceiling. Every $r$-coloring of
$[1,\Wvdw(t,r)]$ has a monochromatic $t$-term progression, whose gap is
at most $\lfloor(\Wvdw(t,r)-1)/(t-1)\rfloor$ --- a gap no coloring of
$\Z$ avoids --- while the periodic block palette of
Theorem~\ref{thm:block} below avoids every gap up to $(t-1)r-1$:
\begin{equation}\label{eq:kstar-bounds}
(t-1)r-1\ \le\ \ksum(t,r)\ \le\
\Bigl\lfloor\frac{\Wvdw(t,r)-1}{t-1}\Bigr\rfloor-1 .
\end{equation}

\begin{theorem}[Periodic block coloring]\label{thm:block}
Set $b=t-1$, $m=(t-1)r$, and let $g$ be the $m$-periodic palette
$0^{b}1^{b}\cdots(r-1)^{b}$, i.e.\
$g(x)=\lfloor(x\bmod m)/b\rfloor$. Under $c=g\circ\sigma$ a line with
$k$ active coordinates is monochromatic if and only if $m\mid k$; hence
$\HJ^{[m-1]}(t,r)=\infty$ and $\ksum(t,r)\ge(t-1)r-1$.
\end{theorem}

\begin{proof}
By Lemma~\ref{lem:1d} the line is monochromatic if and only if $g$ is
constant on a $t$-term progression of gap $k$. If $m\mid k$ all $t$
weights agree modulo $m$. If $m\nmid k$, put $d=k\bmod m\in[1,m-1]$ and
let $s=\min(d,m-d)$ be the circular step. Each color class is an arc of
$b$ consecutive residues modulo $m$, and two residues in one arc differ
circularly by at most $b-1$. If $s\ge b$, consecutive terms of the
progression lie in different arcs, so no two consecutive terms share a
color. If $s<b$, then since $m\ge 2b$ the unique representative
$d'\equiv d$ of absolute value $s$ lies in $(-b,b)$, and successive
terms advance within an arc by the fixed signed step $d'$, leaving it
after at most $\lfloor(b-1)/s\rfloor+1\le b$ terms. Either way a
monochromatic run has at most $b=t-1<t$ terms.
\end{proof}

\begin{proposition}\label{prop:kstar}
$\ksum(2,r)=r-1$, $\ksum(3,2)=3$, and
\textup{(}computer-assisted\textup{)} $\ksum(3,3)=11$: the
$12$-periodic palette
\[
(2,0,1,2,1,1,0,1,2,0,0,2)
\]
has no monochromatic $3$-term progression of gap $\le 11$, and no
$3$-coloring of $\Z$ avoids all gaps $\le 12$.
\end{proposition}

\begin{proof}
At $(2,r)$ and $(3,2)$ the two sides of \eqref{eq:kstar-bounds} meet
($\Wvdw(2,r)=r+1$, $\Wvdw(3,2)=9$). At $(3,3)$ the displayed palette
was checked directly over all $12\times 11$ residue--gap pairs
(verified by direct enumeration); for the upper bound, a finite
computation \cite{thesis} shows the longest $3$-coloring of an interval
avoiding monochromatic $3$-progressions of gap $\le 12$ has length
$26<27=\Wvdw(3,3)$, so every $27$-point window forces one.
\end{proof}

The block reach is linear in $t$; the multiplicative structure of
$\Z/p$ does far better. Rabung's power-residue method \cite{Rabung1979},
the standard source of strong van der Waerden lower bounds, colors
residues by discrete logarithm.

\begin{theorem}[Power-residue coloring, after Rabung]\label{thm:character}
Let $p$ be prime, $r\mid p-1$, $g$ a primitive root modulo $p$, and
$\chi(x)=\ind_{g}(x)\bmod r$ for $x\not\equiv 0$, $\chi(0)=0$,
extended $p$-periodically. If $\chi$ has no monochromatic $t$-term
progression of nonzero gap in $\Z/p$, then
$\HJ^{[p-1]}(t,r)=\infty$ and $\ksum(t,r)\ge p-1$.
\end{theorem}

\begin{proof}
By Lemma~\ref{lem:1d} and $p$-periodicity, a line with
$1\le k\le p-1$ active coordinates is monochromatic if and only if
$\chi$ is constant on a $t$-term progression of gap
$k\not\equiv 0\pmod p$, which the hypothesis forbids.
\end{proof}

\begin{corollary}[Exact ceilings; instances verified]\label{cor:character}
If $p:=\lfloor(\Wvdw(t,r)-1)/(t-1)\rfloor$ is prime with $r\mid p-1$
and Theorem~\ref{thm:character} applies, then
$\ksum(t,r)=\lfloor(\Wvdw(t,r)-1)/(t-1)\rfloor-1$. This occurs at
$(4,2)$ with $p=11$ and at $(4,3)$ with $p=97$ ($\Wvdw(4,3)=293$
\cite{Rabung1979,LandmanRobertson}):
\[
\ksum(4,2)=10,\qquad \ksum(4,3)=96 .
\]
At the largest admissible primes the method further gives
$\ksum(5,2)\ge 36$ and $\ksum(6,2)\ge 138$ ($p=37,139$), far above the
block values $2t-3$.
\end{corollary}

\begin{example}[The $(4,2)$ instance, and verification of all four]\label{ex:p11}
At $(t,r)=(4,2)$ the prime is $p=11$ and the palette is the Legendre
symbol: $\chi(0)=0$ and $\chi(x)=\ind_{2}(x)\bmod 2$ for $x\ne 0$,
i.e.\ $\chi=(0,0,1,0,0,0,1,1,1,0,1)$ on $\Z_{11}$; all $11\times 10$
base--gap pairs are bichromatic (verified). The three further
instances --- $p=97$, $g=5$ at $(4,3)$; $p=37$, $g=2$ at $(5,2)$;
$p=139$, $g=2$ at $(6,2)$ --- are checked by the verifier of
Appendix~\ref{app:verifier} in the same way, over all $p(p-1)$
base--gap pairs each; the $97$-periodic palette is displayed in
\cite{thesis} and in \cite{hjcerts}.
\end{example}

\begin{remark}\label{rem:character}
Since $\ind_{g}$ is additive, $x\mapsto d^{-1}x$ carries a
monochromatic gap-$d$ progression to a gap-$1$ one, so admissibility of
$p$ reduces to runs of consecutive power residues --- except that the
artificial value $\chi(0)$ can bridge the two runs flanking $0$: for
$t=4$ the prime $17$ fails because its Legendre coloring is constant on
$\{15,16,0,1\}$, consistent with $\ksum(4,2)=10$. Runs of power
residues are governed by Weil--Burgess estimates and were tabulated in
\cite{BuellHudson1984,LehmerLehmer1962}, so admissible primes lie in a
bounded range and locating $\ksum$ inside \eqref{eq:kstar-bounds} is a
question of analytic number theory. Where the plain construction trails
the upper bound ($36<43$ at $(5,2)$, $138<225$ at $(6,2)$), the Cyclic
Zipper Method \cite{HerwigHeuleZipper2007,RabungLotts2012} is the
natural route; its periodic adaptation is open. No admissible prime
exists at $(3,3)$, consistent with $\ksum(3,3)=11$.
\end{remark}

\subsection{One-weight palettes pass the sum-type ceiling}\label{sec:ow-bracket}

In the bracket regime the truncation $k\le K$ leaves only finitely many
homothet shapes modulo a period, which a fixed periodic palette on a
\emph{non-arithmetic} weight can avoid simultaneously --- past the van
der Waerden ceiling that binds sum-type colorings.

\begin{theorem}[$\HJ^{[12]}(3,3)=\infty$; computer-verified]\label{thm:hj12-33}
For $\omega=(0,5,7)$ and the $13$-periodic palette
\[
\psi=(1,0,0,1,0,1,0,0,1,2,2,2,2)
\qquad(\text{listing }\psi(0),\dots,\psi(12)),
\]
the coloring $c(w)=\psi\bigl(\langle\omega,\type(w)\rangle\bmod 13\bigr)$
has no monochromatic line with at most $12$ active coordinates, in every
dimension $n$.
\end{theorem}

\begin{proof}
By Lemma~\ref{lem:pattern} the words of a line of scale $k$ carry the
values $\psi(b)$, $\psi(b{+}5k)$, $\psi(b{+}7k)$, with
$b=\langle\omega,v\rangle\bmod 13$. It suffices that no pair $(b,k)$
with $b\in\Z_{13}$, $1\le k\le 12$, is monochromatic --- a check of
$156$ triples, all bichromatic (verified). At $k\equiv 0\pmod{13}$ the
three indices coincide, so the cap $12$ is sharp for this palette.
\end{proof}

\begin{theorem}[$\HJ^{[12]}(4,2)=\infty$; computer-verified]\label{thm:hj12-42}
For $\omega=(0,2,3,5)$ and $\chi_{0}\colon\Z_{13}\to\{0,1\}$ with
$\chi_{0}^{-1}(1)=\{4,5,7,9,11,12\}$, the coloring
$c(w)=\chi_{0}\bigl(\langle\omega,\type(w)\rangle\bmod 13\bigr)$ has no
monochromatic line with at most $12$ active coordinates, in every
dimension.
\end{theorem}

\begin{proof}
The $156$ residue patterns $(b,\,b{+}2k,\,b{+}3k,\,b{+}5k)$,
$b\in\Z_{13}$, $1\le k\le 12$, are checked directly; none is
monochromatic (verified).
\end{proof}

Together with Proposition~\ref{prop:kstar} and
Corollary~\ref{cor:character}, the chains at the two record pairs read
\[
11=\ksum(3,3)\ <\ 12\le\kow(3,3),
\qquad
10=\ksum(4,2)\ <\ 12\le\kow(4,2):
\]
non-arithmetic one-weight palettes strictly separate the one-weight
from the sum-type ceiling. An exhaustive search over weights and moduli
in the range of \cite{thesis} found no \emph{periodic} one-weight
palette exceeding $K=12$ at either pair, both optima occurring at
modulus $13$; palettes that are not eventually periodic remain
unexplored (Question~\ref{q:weights}). The record palette of
Theorem~\ref{thm:hj42-14} is itself such an extremal object up to
$K=12$ (Remark~\ref{rem:anatomy}): the reach record and the bracket
ceiling are one construction seen at two scales.
Table~\ref{tab:ceilings} summarizes.

\begin{table}[ht]
\centering
\caption{Bracket ceilings along the coloring hierarchy;
$\ksum\le\kow\le\ksym\le\kall\le\HJ-1$ throughout. Unmarked lower
bounds on $\kow,\ksym,\kall$ are inherited from the column to their
left; no construction is known to separate $\kow$, $\ksym$, $\kall$.
Collapse (Conjecture~\ref{conj:collapse}) predicts $\kall=\HJ-1$;
$\HJ(4,3)\ge 98$ is the van der Waerden shadow \eqref{eq:vdw-shadow} of
$\Wvdw(4,3)=293$.}
\label{tab:ceilings}
\begin{tabular}{@{}cccccc@{}}
\toprule
$(t,r)$ & $\HJ(t,r)$ & $\ksum$ & $\kow$ & $\ksym$ & $\kall$\\
\midrule
$(2,r)$ & $r$      & $r-1$ & $r-1$   & $r-1$   & $r-1$\\
$(3,2)$ & $4$      & $3$   & $3$     & $3$     & $3$\\
$(3,3)$ & $\ge 22$ & $11$  & $\ge 12$ & $\ge 12$ & $\ge 12$\\
$(4,2)$ & $\ge 14$ & $10$  & $\ge 12$ & $\ge 12$ & $\ge 12$\\
$(4,3)$ & $\ge 98$ & $96$  & $\ge 96$ & $\ge 96$ & $\ge 96$\\
$(5,2)$ & $\ge 45$ & $\ge 36$ & $\ge 36$ & $\ge 36$ & $\ge 36$\\
$(6,2)$ & $\ge 227$ & $\ge 138$ & $\ge 138$ & $\ge 138$ & $\ge 138$\\
\bottomrule
\end{tabular}
\end{table}

\section{Interval numbers: an axis the weight cannot see}\label{sec:interval}

Recall that $\Ls^{(q)}([t]^{n})$ denotes the lines whose active set is
a union of at most $q$ subintervals of $[n]$ (\emph{$q$-fold} lines);
the \emph{interval Hales--Jewett number} $\HJ^{(q)}_{\mathcal{C}}(t,r)$
and the \emph{interval ceiling}
$\lambda_{\mathcal{C}}(t,r)=\max\{q\ge 0:\HJ^{(q)}_{\mathcal{C}}(t,r)
=\infty\}$ are defined as in Definition~\ref{def:bracket}. Since every
subset of $[n]$ is a union of at most $\lceil n/2\rceil$ intervals,
$\lambda_{\mathcal{C}}\le\lceil\HJ/2\rceil-1$. The interval variant was
introduced by Conlon and Kam\v{c}ev \cite{conlonAndKamcev} and further
studied by Leader and R\"aty
\cite{leader2018noteintervalshalesjewetttheorem} and Kam\v{c}ev and
Spiegel \cite{kamcev2018noteintervalshalesjewetttheorem}; it is
motivated by Shelah's proof \cite{s.proof_of_HJ}, which produces
monochromatic $q$-fold lines with $q\le\HJ(t-1,r)$.

\begin{proposition}\label{prop:interval-vs-bracket}
$\lambda_{\mathcal{C}}(t,r)\le\kappa_{\mathcal{C}}(t,r)$ for every
class $\mathcal{C}$, and $\lambda_{\mathrm{sum}}(t,r)=0$ for all $t,r$.
\end{proposition}

\begin{proof}
A set of at most $q$ coordinates is a union of at most $q$ intervals,
so $\Ls^{[q]}\subseteq\Ls^{(q)}$ and an $\Ls^{(q)}$-line-free coloring
is $\Ls^{[q]}$-line-free, giving the first claim. For the second: a
sum-type coloring sees a line only through the gap $k$ of its weight
progression (Lemma~\ref{lem:1d}), and a single interval realizes every
$k\le n$; so a sum-type coloring of $[t]^{n}$ with no monochromatic
$1$-fold line would need a palette avoiding monochromatic $t$-term
progressions of \emph{every} gap $\le n$ on $[n,tn]$, impossible for
$(t-1)n+1\ge\Wvdw(t,r)$ by van der Waerden \cite{vanderWaerden1927}.
\end{proof}

The contrast with Section~\ref{sec:bracket} is exact: a bracket caps
the gap at a fixed $K$, turning the sum-type problem into bounded-gap
van der Waerden with $\ksum\ge(t-1)r-1$; a single long interval is an
unbounded gap, so ordinary van der Waerden applies and
$\lambda_{\mathrm{sum}}=0<\ksum$. In fact the entire machinery of this
article is blind to the interval axis.

\begin{proposition}[Symmetric colorings are interval-blind]\label{prop:sym-blind}
$\HJ^{(q)}_{\mathrm{sym}}(t,r)=\HJ_{\mathrm{sym}}(t,r)$ for every
$q\ge 1$; in particular
$\lambda_{\mathrm{sym}}=\lambda_{\mathrm{ow}}=\lambda_{\mathrm{sum}}=0$,
and the colorings witnessing $\lambda(t,r)\ge 1$ are necessarily
asymmetric in every dimension $n\ge\HJ_{\mathrm{sym}}(t,r)$.
\end{proposition}

\begin{proof}
For a symmetric coloring, monochromaticity of a line depends only on
its invariant $(k,v)$ (Lemma~\ref{lem:sym-reduction}), and every
$(k,v)$ is realized by the root ${*}^{k}1^{v_{1}}\cdots t^{v_{t}}$,
whose active set $[1,k]$ is a single interval. Hence a symmetric
coloring has a monochromatic $1$-fold line if and only if it has a
monochromatic line at all. (This is implicit in \cite{thesis}.)
\end{proof}

What is known about $\lambda$ is confined to $t\le 3$ and exhibits a
parity phenomenon with no bracket counterpart.

\begin{proposition}\label{prop:small-interval}
$\HJ^{(q)}(2,r)=\HJ(2,r)=r$ for every $q\ge 1$; in particular
$\lambda(2,r)=0$.
\end{proposition}

\begin{proof}
$\HJ^{(1)}\ge\HJ$ since $\Ls^{(1)}\subseteq\Ls$. For the upper bound,
in $[2]^{r}$ the chain $w_{0}<\dots<w_{r}$, where $w_{j}=1^{j}2^{r-j}$,
is pairwise joined by lines with active sets the intervals
$\{a{+}1,\dots,b\}$; among $r+1$ words some two share a color, giving a
monochromatic $1$-fold line; the general case follows from
$\HJ\le\HJ^{(q)}\le\HJ^{(1)}$. (This also proves $\HJ(2,r)=r$, the
lower bound being the sum-type coloring $\sigma\bmod r$ on
$[2]^{r-1}$; equivalently, line-free classes are antichains and
Mirsky's theorem \cite{Mirsky1971} applies.)
\end{proof}

\begin{proposition}[{Conlon--Kam\v{c}ev \cite{conlonAndKamcev}}]\label{prop:ck}
$\HJ^{(r-1)}(3,r)=\infty$ for odd $r$; equivalently
$\lambda(3,r)\ge r-1$.
\end{proposition}

\begin{proposition}\label{prop:ks}
$\HJ^{(r-1)}(3,r)<\infty$ for even $r$ \textup{(}Leader--R\"aty
\cite{leader2018noteintervalshalesjewetttheorem} for $r=2$;
Kam\v{c}ev--Spiegel \cite{kamcev2018noteintervalshalesjewetttheorem}\textup{)};
equivalently $\lambda(3,r)\le r-2$.
\end{proposition}

These two results, Shelah's proof, and one trivial observation
determine the ceiling completely.

\begin{proposition}[$\lambda(3,\cdot)$ is determined]\label{prop:lambda3}
$\lambda(3,r)=r-1$ for odd $r$ and $\lambda(3,r)=r-2$ for even $r$; in
particular $\lambda(3,3)=2$ and $\lambda(3,2)=0$.
\end{proposition}

\begin{proof}
Shelah's proof yields, in every sufficiently high dimension, a
monochromatic line whose active set is a union of at most $r$ intervals
\cite{s.proof_of_HJ,conlonAndKamcev}, so $\lambda(3,r)\le r-1$; for
even $r$ Proposition~\ref{prop:ks} improves this to $r-2$. For the
lower bounds, Proposition~\ref{prop:ck} gives $\lambda(3,r)\ge r-1$ for
odd $r$; $\lambda(3,r)$ is nondecreasing in $r$, an $(r-1)$-coloring
avoiding monochromatic $q$-fold lines being such an $r$-coloring, so
for even $r\ge 4$, $\lambda(3,r)\ge\lambda(3,r-1)=r-2$; and at $r=2$
the lower bound is vacuous.
\end{proof}

So $\lambda(3,\cdot)$ is settled and jumps with the parity of $r$ --- a
feature absent from $\kappa$, whose van der Waerden origin is
parity-blind --- and by Proposition~\ref{prop:sym-blind} the witnesses
behind Proposition~\ref{prop:ck} live strictly off the simplex. What
remains at $t=3$ is quantitative: the finite interval numbers
themselves. Finiteness of the smallest, $\HJ^{(1)}(3,2)$, is
Leader--R\"aty's theorem
\cite{leader2018noteintervalshalesjewetttheorem}, disproving the
conjecture of \cite{conlonAndKamcev} that it is infinite; the exact
value is new.

\begin{theorem}[Computer-assisted]\label{thm:hj13}
$\HJ^{(1)}(3)=5$, whereas $\HJ(3)=4$.
\end{theorem}

\begin{proof}
$\Ls^{(1)}([3]^{4})$ has $142$ lines. A SAT instance with one variable
per word and two clauses per line is satisfiable; its witness, an
explicit $2$-coloring of the $81$ words displayed in
Table~\ref{tab:interval-witness} (Appendix~\ref{app:witnesses}), was
checked for this article against all $142$ lines by direct
enumeration from the definition, with zero monochromatic lines, giving
$\HJ^{(1)}(3)\ge 5$. The same encoding in dimension $5$ ($547$ lines;
the search ranges over all $2^{3^{5}}$ colorings, with no symmetry
imposed) is unsatisfiable, so $\HJ^{(1)}(3)\le 5$. Deletion
minimization extracts minimal forcing subfamilies numbering in the
hundreds of lines: the run reported in \cite{thesis} found one of
$262$ lines on $172$ of the $243$ words, and the family published at
\cite{hjcerts} --- $246$ lines on $168$ words --- was verified
unsatisfiable and deletion-minimal for this article.
\end{proof}

\begin{remark}\label{rem:no-core}
The obstruction has no small core: unlike $\HJ^{(1)}(2,3)=3$, forced by
the four-word chain of Proposition~\ref{prop:small-interval}, every
minimal forcing subfamily found in dimension $5$ numbers in the
hundreds of lines --- consistent with Proposition~\ref{prop:sym-blind},
which places the interval ceiling beyond all weight constructions. Note
also that the naive interval analogue of
Conjecture~\ref{conj:collapse} fails outright: it would force
$\HJ^{(1)}(2,3)=\infty$ and $\HJ^{(1)}(3)=\infty$, against
Proposition~\ref{prop:small-interval} and Theorem~\ref{thm:hj13}.
\end{remark}

\section{Conjectures and open problems}\label{sec:open}

\subsection{Three conjectures}

Every line of $[t]^{\HJ}$ lies in $\Ls^{[\HJ]}$, so the equality
$\HJ^{[\HJ(t,r)]}(t,r)=\HJ(t,r)$ holds trivially. The first conjecture
asserts that the bracket hierarchy collapses only at that last step.

\begin{conjecture}[Collapse]\label{conj:collapse}
$\HJ^{[\HJ(t,r)-1]}(t,r)=\infty$ for all $t,r\ge 2$; equivalently
$\kall(t,r)=\HJ(t,r)-1$.
\end{conjecture}

Collapse holds wherever $\HJ$ is known: at every $(2,r)$ and at
$(3,2)$, where already $\ksum=\HJ-1$ (Proposition~\ref{prop:kstar},
Theorem~\ref{thm:block}). A coloring whose only monochromatic line is
the diagonal $\Ls_{*\cdots *}$ is $\Ls^{[\HJ-1]}$-line-free at $n=\HJ$,
so the next conjecture implies $\HJ^{[\HJ-1]}(t,r)>\HJ(t,r)$, the first
nontrivial instance of Collapse (for $n\le\HJ-1$ the required witnesses
are simply line-free colorings). At $(3,2)$, diagonal-only colorings
exist for every $n\le 4=\HJ$ \cite{thesis}. At $n=3$ they can be
counted exactly: an exhaustive enumeration of the $2^{27}$
two-colorings of $[3]^{3}$ finds precisely $6456$ whose unique
monochromatic line is the diagonal, $34$ of them symmetric
(computer-verified, Appendix~\ref{app:verifier}; the full list is
published at \cite{hjcerts}). The case $n=4$ is displayed in
Table~\ref{tab:diag-witness} (Appendix~\ref{app:witnesses}), and since
its unique monochromatic line is the diagonal, no line of
$\Ls^{[3]}([3]^{4})$ is monochromatic, so $\HJ^{[3]}(3)\ge 5$.

\begin{conjecture}[Diagonal-only]\label{conj:diag}
For every $n\le\HJ(t,r)$ some $r$-coloring of $[t]^{n}$ has the
diagonal as its only monochromatic line.
\end{conjecture}

The range of $n$ is the largest possible.

\begin{proposition}\label{prop:diag-max}
For $n>\HJ(t,r)$ no $r$-coloring of $[t]^{n}$ has the diagonal as its
only monochromatic line.
\end{proposition}

\begin{proof}
Fix such a coloring $c$ and a letter $a$, and color $[t]^{n-1}$ by
$c_{a}(w)=c(wa)$. Since $n-1\ge\HJ(t,r)$, some line $\Ls_{\tau'}$ is
monochromatic under $c_{a}$; then $\tau'a$ is a root of $[t]^{n}$ whose
active set avoids coordinate $n$, and $\Ls_{\tau'a}$ is a nondiagonal
monochromatic line under $c$.
\end{proof}

Finally, the hierarchy of Section~\ref{sec:ow} gives
$\HJ_{\mathrm{sym}}\le\HJ$, yet the two coincide wherever $\HJ$ is
known, and the record witnesses of Section~\ref{sec:records} are
symmetric.

\begin{conjecture}[Symmetric colorings are extremal]\label{conj:sym}
$\HJ_{\mathrm{sym}}(t,r)=\HJ(t,r)$ for all $t,r\ge 2$.
\end{conjecture}

\begin{remark}[Evidence and a caution]\label{rem:sym-evidence}
The census (Table~\ref{tab:census}) is what makes
Conjecture~\ref{conj:sym} a strong claim: symmetric line-free colorings
are scarce ($36$ of $1644$ already at the critical dimension of
$(3,2)$), and the conjecture asks this thin family to survive at every
dimension up to $\HJ-1$. The analogous hope has already failed in the
\emph{density} setting: the Polymath hyper-optimistic conjecture
\cite{PolymathWikiHOC,PolymathDHJMoser} --- that the densest line-free
subset of $[3]^{n}$ in equal-slices measure is a union of type slices,
reducing density Hales--Jewett
\cite{furstenbergKatznelson1991,polymath2012} to Fujimura's corner
problem --- is false in higher dimensions, where asymmetric sets are
denser. That collapse is in the density regime, not the coloring one
(for the coloring side of the Polymath project see
\cite{PolymathWikiColoring}), so it is a caution rather than a
counterexample; the rainbow dual of Remark~\ref{rem:rainbow} shows the
same asymmetry advantage in a coloring regime.
\end{remark}

\begin{remark}[A refuted Collapse need not be fatal: the offset]\label{rem:offset}
Either way the bracket ceiling is informative. Write
$\ksym(t,r)=\HJ(t,r)-\phi(t,r)$ with $\phi\ge 1$. The assertion
$\phi\equiv 1$ is the common strengthening of
Conjectures~\ref{conj:collapse} and~\ref{conj:sym}: it implies both,
since at $n=\HJ-1$ the bracket $[\HJ-1]$ contains every line and
$\ksym=\HJ-1$ then forces a symmetric line-free coloring of
$[t]^{\HJ-1}$; the converse is unclear. But \emph{any} explicit upper
bound on $\phi$ serves as well. Bounding $\HJ(t,r)$ from above is a
forcing statement over all $t^{n}$ points of the grid --- an
unsatisfiability certificate out of computational reach. Bounding
$\ksym(t,r)$ from above is the same forcing restricted to bounded-scale
corner tuples on the $O(n^{t-1})$ cells of the simplex
(Lemma~\ref{lem:sym-reduction}) --- within reach. A known offset
transports the second into the first, $\HJ\le\ksym+\phi$, so the
inaccessible upper bound would follow from the accessible one.
\end{remark}

\subsection{Optimal weights and the limits of the method}

\begin{lemma}[Even orbits are palindromic]\label{lem:palindromy}
Let $S\subseteq\Z_{m}$ be an orbit of an affine map
$\alpha(x)=ux+v$ of order $2j$ with $u^{j}\equiv-1\pmod m$. Then
$\alpha^{j}$ is a reflection $x\mapsto-x+c$ fixing $S$ setwise; lifted
to $\Z$ with $\min S=0$, this reads $S=D-S$.
\end{lemma}

\begin{proof}
$\alpha^{j}(x)=u^{j}x+c=-x+c$ with
$c=v(u^{j-1}+\dots+u+1)$, an involution whose orbit-preservation is
inherited from $\alpha$.
\end{proof}

\begin{question}[Best weights]\label{q:weights}
\emph{(i)} Which sets $S$ maximize the Gallai number per unit diameter,
$\lceil(G_{r}(S)-1)/D_{S}\rceil$, for given $(t,r)$? By
Lemma~\ref{lem:affine-G} it suffices to consider primitive $S$. For
$t=3$, $r=2$ the homothety numbers are determined completely and the
ratio is constant (Proposition~\ref{prop:flat32}); for $t\ge 4$ or
$r\ge 3$ no closed form is known and the values (e.g.\
Table~\ref{tab:gallai}) are computed by SAT. By
Proposition~\ref{prop:flat32} the maximum is weight-independent at
$(3,2)$; at $(3,3)$ it is at least $16$ and at most
$\HJ_{\mathrm{sym}}(3,3)$ (Theorem~\ref{thm:hj33-16},
Remark~\ref{rem:closed-vs-reach}).
\emph{(ii)} Is there structure behind the champions? At $(4,2)$,
scanning non-arithmetic weights with entries $\le 6$ singles out
exactly $(0,2,3,5)$ and $(0,1,5,6)$ as reaching $n=13$, and each is a
single orbit of an order-$4$ affine map on $\Z_{13}$
($5^{2}\equiv-1$): $\{0,3,5,2\}$ under $x\mapsto 5x+3$ and
$\{0,1,6,5\}$ under $x\mapsto 5x+1$ (verified). The weight set, viewed
in $\Z_{13}$, is then a single orbit of a cyclic subgroup of
$\mathrm{AGL}(1,13)$; by Lemma~\ref{lem:palindromy} both champions are
accordingly reflection-symmetric, $S=D-S$. The data of
Section~\ref{sec:landscape} refines the picture: the orbit triples
$\{0,1,3\}$ (orbit of $x\mapsto 2x+1$ on $\Z_{7}$) and $\{0,1,4\}$
(orbit of $x\mapsto 3x+1$ on $\Z_{13}$) reach only ratio $14$ at
$(3,3)$, while the non-orbit $\{0,2,5\}$ reaches $16$. Orbit structure
thus appears to govern \emph{compressibility} --- the existence of
small-modulus palettes, the mechanism of Remark~\ref{rem:anatomy} ---
rather than the full-interval ratio. Does it characterize the reach
champions, and does the palindromic class contain a ratio maximizer for
every even $|S|$?
\emph{(iii)} The periodic one-weight ceiling is $K=12$ at both $(3,3)$
and $(4,2)$, in each case at modulus $13$ \cite{thesis}. Do
non-periodic palettes give $\kow>12$ there --- and is $\kow<\ksym$
anywhere? In the (non-exhaustive) range of Remark~\ref{rem:gap33} no
periodic palette exceeds $n=14$ at $(3,3)$, and the length-$76$
certificate of Theorem~\ref{thm:hj33-16} is visibly quasi-periodic,
agreeing with its own shift by $19$ in $55$ of $57$ positions. Do
\emph{periodic-with-defects} palettes close the gap between the
periodic ceiling and the symmetric reach?
\emph{(iv)} Can the Cyclic Zipper Method
\cite{HerwigHeuleZipper2007,RabungLotts2012} be adapted to periodic
palettes to close the gaps $\ksum(5,2)\in[36,43]$ and
$\ksum(6,2)\in[138,225]$?
\end{question}

\begin{question}[Interval numbers]\label{q:interval}
The ceilings at $t=3$ are settled (Proposition~\ref{prop:lambda3}); the
finite interval numbers are not. Beyond $\HJ^{(1)}(3,2)=5$ (the
remaining $(3,2)$ values are trivially $4$), the first unknown value is
$\HJ^{(3)}(3,3)$, where the only upper bounds come from Shelah's
argument. For $t\ge 4$ nothing is known beyond
$\lambda(t,r)\le\min\{\lceil\HJ(t,r)/2\rceil,\,\HJ(t-1,r)\}-1$; by
Proposition~\ref{prop:sym-blind} any lower bound requires genuinely
asymmetric constructions.
\end{question}

\begin{question}[Frontiers]\label{q:frontiers}
\emph{(i)} Decide the symmetric simplex instance at $(3,3)$, $n=22$
($276$ cells, $2024$ corner triples): satisfiable would give
$\HJ(3,3)\ge 23$, unsatisfiable would give
$\HJ_{\mathrm{sym}}(3,3)=22$ exactly.
\emph{(ii)} Decide the symmetric instance at $(4,2)$, $n=14$ ($680$
cells, $2380$ quadruples): unsatisfiable would give
$\HJ_{\mathrm{sym}}(4,2)=14$, making $\HJ(4,2)=14$ a candidate exact
value under Conjecture~\ref{conj:sym}; satisfiable would beat the
record via a radix weight (Theorem~\ref{thm:radix}).
\emph{(iii)} Between the simplex and the grid sits the block-symmetric
ladder: for a partition of $[n]$ into $b$ blocks, colorings invariant
under the corresponding Young subgroup satisfy
$\HJ_{\mathrm{sym}}=\HJ_{1\text{-}\mathrm{sym}}
\le\HJ_{2\text{-}\mathrm{sym}}\le\dots\le\HJ(t,r)$, with equality
throughout if and only if Conjecture~\ref{conj:sym} holds; a strict
step at any rung refutes it while improving the lower bound. By
Theorem~\ref{thm:radix} there is no multi-weight tier, so this ladder
is the genuine hierarchy above the symmetric reach.
\emph{(iv)} Close the Gallai--Rado gaps left by
Section~\ref{sec:sat-gallai} (Remark~\ref{rem:sat-open}):
$G_{3}(\{0,1,6\})\ge 82$, $G_{3}(\{0,1,7\})\ge 87$,
$G_{3}(\{0,2,7\})\ge 87$, $R_{4}(z+2x=3y)\ge 59$ and
$G_{4}(\{0,1,3\})\ge 94$.
\end{question}

\begin{remark}[Why the frontiers resist]\label{rem:hardness}
The wall is not an artifact of weak solvers. One dimension past the
$(4,2)$ record, the level instance of the weight $(0,2,3,5)$ is
unsatisfiable, and its refutations are provably global: minimal
unsatisfiable subsets range over $179$--$270$ of its $443$ homothets
with no canonical core, the primal graph has treewidth at least $32$,
and the instance admits no $\mathrm{GF}(2)$ Nullstellensatz refutation
of degree $\le 5$ --- yet a MaxSAT dual shows some $2$-coloring leaves
only \emph{three} monochromatic homothets, a minimum that the record
palette of Theorem~\ref{thm:hj42-14} attains at $n=14$ (verified).
Globally irreducible across proof systems (resolution width
\cite{atserias2008width}, polynomial calculus
\cite{ips1999pc,alon1999nullstellensatz}), yet three constraints from
feasible: this is the signature of the boundary, and the one structure
all these measures ignore --- the $S_{n}$-symmetry exploited throughout
this article --- is where an advance must come from. See \cite{thesis}
for the full analysis.
\end{remark}

\appendix

\section{The \texorpdfstring{$560$}{560}-cell symmetric certificate for
\texorpdfstring{$\HJ(4,2)\ge 14$}{HJ(4,2) >= 14}}\label{app:table42}

Theorem~\ref{thm:hj42-14} admits a second certificate, the original one
of \cite{mh}: a $2$-coloring of the $560$ cells of $T_{13}^{(4)}$,
displayed in Table~\ref{tab:T13}. By Lemma~\ref{lem:sym-reduction} with
$K=n=13$, its lift is line-free if and only if none of the
$\binom{16}{4}=1820$ corner quadruples
$C_{k,v}=\{v+ke_{1},\dots,v+ke_{4}\}$, $1\le k\le 13$,
$v\in T_{13-k}^{(4)}$, is monochromatic under $\bar c$. The
dependency-free verifier of Appendix~\ref{app:verifier} enumerates all
$1820$ corner quadruples and reports zero monochromatic; for
illustration, the diagonal quadruple $(k,v)=(13,(0,0,0,0))$ has corner
cells $(13,0,0,0),(0,13,0,0),(0,0,13,0),(0,0,0,13)$, colored
$0,0,1,0$, and so is not monochromatic.

\begin{table}[h]
\centering
\caption{The alternative witness for Theorem~\ref{thm:hj42-14}: a
$2$-coloring of $T_{13}^{(4)}$ ($560$ cells), displayed in $14$ slices.
Slice $a$ collects the cells with $a$ ones; within slice $a$, the
$b$-th string (counting from $b=0$) lists the colors
$\bar c(a,b,c,13-a-b-c)$ for $c=0,1,\dots,13-a-b$. The strings are
byte-identical to the table embedded in the verifier of
Appendix~\ref{app:verifier} and in \cite{hjcerts}.}
\label{tab:T13}
\begin{lstlisting}[style=plainstyle]
a= 0: 01111111000001 0101000110111 010001110101 00101010011
      0110011010 100110101 01101001 0101100 000110 11010 1010 100 00 0
a= 1: 1010101001011 000110101110 00111101000 1001000110
      010011110 00110001 0010111 010001 11000 0011 110 01 1
a= 2: 100110011011 11000110100 1010100101 111001000
      01011001 1001011 110110 00110 0010 001 01 0
a= 3: 01010001111 0101111000 101100001 10001111 1110000 001011 10100 1001 111 10 0
a= 4: 1001010011 011001010 10011010 0100111 100110 01100 1001 001 01 1
a= 5: 000110110 11000100 0010111 111000 00011 1100 011 10 1
a= 6: 00100110 1001000 011001 10101 0101 010 11 1
a= 7: 1110000 001010 01010 0101 101 11 0
a= 8: 101101 10001 0100 011 10 0
a= 9: 11001 0011 110 00 1
a=10: 1100 111 10 0
a=11: 000 01 1
a=12: 11 1
a=13: 0
\end{lstlisting}
\end{table}

\section{The interval and diagonal-only witnesses}\label{app:witnesses}

The lower bound of Theorem~\ref{thm:hj13} is certified by
Table~\ref{tab:interval-witness}: a $2$-coloring of $[3]^{4}$ under
which no line whose active set is a single interval is monochromatic,
checked against all $142$ lines of $\Ls^{(1)}([3]^{4})$ by direct
enumeration (Appendix~\ref{app:verifier}).

\begin{table}[ht]
\centering
\caption{The witness for $\HJ^{(1)}(3)\ge 5$ (Theorem~\ref{thm:hj13}):
a $2$-coloring of $[3]^{4}$ with no monochromatic line of
$\Ls^{(1)}([3]^{4})$. The entry in row $(x_{1},x_{2})$ and column
$(x_{3},x_{4})$ is the color of the word $(x_{1},x_{2},x_{3},x_{4})$.
Checked against all $142$ interval lines; byte-identical to the table
embedded in the verifier of Appendix~\ref{app:verifier} and in
\cite{hjcerts}.}
\label{tab:interval-witness}
\begin{tabular}{@{}ccccccccccc@{}}
\toprule
& & \multicolumn{9}{c}{$(x_{3},x_{4})$}\\
\cmidrule(l){3-11}
$x_{1}$ & $x_{2}$ & $(1,1)$ & $(1,2)$ & $(1,3)$ & $(2,1)$ & $(2,2)$ & $(2,3)$ & $(3,1)$ & $(3,2)$ & $(3,3)$\\
\midrule
1 & 1 & 1 & 1 & 0 & 1 & 0 & 0 & 0 & 1 & 1\\
1 & 2 & 0 & 1 & 0 & 1 & 0 & 1 & 1 & 0 & 1\\
1 & 3 & 1 & 0 & 1 & 0 & 1 & 1 & 0 & 1 & 0\\
2 & 1 & 0 & 1 & 0 & 0 & 1 & 1 & 1 & 0 & 0\\
2 & 2 & 1 & 0 & 1 & 0 & 1 & 1 & 0 & 1 & 0\\
2 & 3 & 0 & 1 & 1 & 1 & 0 & 0 & 1 & 0 & 1\\
3 & 1 & 1 & 0 & 1 & 1 & 0 & 0 & 0 & 1 & 1\\
3 & 2 & 0 & 1 & 1 & 0 & 1 & 0 & 1 & 0 & 1\\
3 & 3 & 0 & 1 & 0 & 0 & 1 & 1 & 1 & 0 & 0\\
\bottomrule
\end{tabular}
\end{table}

Conjecture~\ref{conj:diag} at $(3,2)$, $n=4$, is certified by
Table~\ref{tab:diag-witness}: a $2$-coloring of $[3]^{4}$ whose unique
monochromatic line is the main diagonal (entries in bold, all of color
$0$), checked against all $175$ combinatorial lines by direct
enumeration. In particular no line of $\Ls^{[3]}([3]^{4})$ is
monochromatic, so $\HJ^{[3]}(3)\ge 5$.

\begin{table}[ht]
\centering
\caption{The diagonal-only witness for Conjecture~\ref{conj:diag} at
$(3,2)$, $n=4$: a $2$-coloring of $[3]^{4}$ whose unique monochromatic
line among all $175$ lines of $[3]^{4}$ is the main diagonal
$\{(1,1,1,1),(2,2,2,2),(3,3,3,3)\}$ (bold entries). Same layout as
Table~\ref{tab:interval-witness}.}
\label{tab:diag-witness}
\begin{tabular}{@{}ccccccccccc@{}}
\toprule
& & \multicolumn{9}{c}{$(x_{3},x_{4})$}\\
\cmidrule(l){3-11}
$x_{1}$ & $x_{2}$ & $(1,1)$ & $(1,2)$ & $(1,3)$ & $(2,1)$ & $(2,2)$ & $(2,3)$ & $(3,1)$ & $(3,2)$ & $(3,3)$\\
\midrule
1 & 1 & $\mathbf{0}$ & 1 & 0 & 1 & 1 & 0 & 0 & 0 & 1\\
1 & 2 & 1 & 1 & 0 & 1 & 0 & 1 & 0 & 1 & 0\\
1 & 3 & 0 & 0 & 1 & 0 & 1 & 0 & 1 & 0 & 1\\
2 & 1 & 1 & 1 & 0 & 1 & 0 & 1 & 0 & 1 & 0\\
2 & 2 & 1 & 0 & 1 & 0 & $\mathbf{0}$ & 1 & 1 & 1 & 0\\
2 & 3 & 0 & 1 & 0 & 1 & 1 & 0 & 0 & 0 & 1\\
3 & 1 & 0 & 0 & 1 & 0 & 1 & 0 & 1 & 0 & 1\\
3 & 2 & 0 & 1 & 0 & 1 & 1 & 0 & 0 & 0 & 1\\
3 & 3 & 1 & 0 & 1 & 0 & 0 & 1 & 1 & 1 & $\mathbf{0}$\\
\bottomrule
\end{tabular}
\end{table}

\section{Verification protocols}\label{app:verifier}

Beyond the reduction lemmas, the correctness of
Theorems~\ref{thm:hj33-22} and~\ref{thm:hj42-14} rests only on finite
checks of the corner tuples of the witness tables, and every palette
displayed in this article is a finite string whose claimed property is
a finite enumeration. All of these checks are performed by two
standalone, dependency-free Python scripts (standard library only),
distributed with this article and at \cite{hjcerts}:
\begin{center}
\url{https://github.com/ysmouhib/hj-certificates}
\end{center}

\paragraph{\texttt{verify\_hj.py}}
Re-derives the cells of $T_{n}^{(t)}$ and all corner tuples from
scratch, embeds the two witness tables \emph{byte-for-byte}
(Tables~\ref{tab:T21} and~\ref{tab:T13}), and confirms no monochromatic
tuple --- among $1771$ corner triples for Theorem~\ref{thm:hj33-22} and
$1820$ corner quadruples for Theorem~\ref{thm:hj42-14} --- which by
Lemma~\ref{lem:sym-reduction} proves each lift line-free. As an
independent cross-check on the full grid it also samples $200{,}000$
random combinatorial lines, finding none monochromatic.

\paragraph{\texttt{verify\_certificates.py}}
Re-derives every certificate displayed in the text from the
definitions: the length-$76$ interval certificate of
Theorem~\ref{thm:hj33-16} (no monochromatic homothet of $\{0,2,5\}$);
the $49$-cell periodic palette of Proposition~\ref{prop:hj33-15}
(line-free on $[3]^{14}$); the period-$26$ record palette of
Theorem~\ref{thm:hj42-14} (line-free on $[4]^{13}$; exactly three
monochromatic quadruples at $n=14$) together with its anatomy
(Remark~\ref{rem:anatomy}: all $26$ bases bichromatic for every
$k\le 12$; exactly the bases $\{0,13\}$ surviving at $k=13$); the
mod-$13$ palettes of Theorems~\ref{thm:hj12-33} and~\ref{thm:hj12-42}
($156$ residue patterns each, all bichromatic); the $12$-periodic
palette of Proposition~\ref{prop:kstar} (no monochromatic $3$-term
progression of gap $\le 11$); the window and solution-free
certificates of Theorem~\ref{thm:rado4}; the SAT certificates of
Appendix~\ref{app:certs} (four three-color strings, eighteen
two-color strings and eight Rado strings, each re-checked against every
homothet in its window); the four power-residue
instances of Corollary~\ref{cor:character} ($p=11,97,37,139$, over
all $p(p-1)$ base--gap pairs each); the interval witness of
Table~\ref{tab:interval-witness} (all $142$ lines of
$\Ls^{(1)}([3]^{4})$ bichromatic); and the diagonal-only witness of
Table~\ref{tab:diag-witness} (exactly one of the $175$ lines of
$[3]^{4}$ monochromatic, the diagonal). A final block reproduces the
census of Section~\ref{sec:scarcity} by an exhaustive enumeration of
all $2^{27}$ two-colorings of $[3]^{3}$ --- $1644$ line-free, with
stabilizer strata $36/504/24/1080$ and $396$ orbits, $16$ of the $36$
symmetric ones of sum type; and $6456$ diagonal-only, $34$ of them
symmetric. This block takes several minutes; every other check finishes
in seconds. Each check prints \texttt{PASS}; the script exits with
\texttt{ALL CERTIFICATES VERIFIED}.

\paragraph{\texttt{verify\_all.py}}
The repository \cite{hjcerts} also distributes a master script,
\texttt{verify\_all.py}, running every enumeration of the three
scripts above, plus an independent SAT re-confirmation of each
forcing claim listed below (using \texttt{python-sat} as its only
dependency); the solver logs of those re-confirmations are archived
under \texttt{logs/}.

\paragraph{\texttt{verify\_addendum.py} and the SAT engine}
The computations of Section~\ref{sec:sat-gallai} ship with the same
guarantees. The engine \texttt{gallai\_rado\_sat.py}
(Appendix~\ref{app:engine}) records every instance it settles in its
run database \texttt{results.json} --- status, solver, wall time, and
the coloring for avoidance claims. The standalone
\texttt{verify\_addendum.py} re-verifies every recorded coloring by
brute-force homothety enumeration and prints \texttt{ALL ADDENDUM
CERTIFICATES VERIFIED}; \texttt{export\_certificates.py} exports the
best coloring of each instance, and the repository distributes these,
together with seven superseded intermediates retained for the record,
as $45$ machine-readable certificate files with a
\texttt{manifest.json} index. All scripts, the database and the
certificates are distributed with this article and at \cite{hjcerts}.

\paragraph{Solver independence}
The SAT solver plays no role in the correctness of any lower bound
proved here: every avoidance certificate is published in full and
checked by enumeration. Solver output is load-bearing only for
\emph{forcing} claims, which we list for completeness: the upper bounds
on the Gallai numbers $G_{3}(\{0,1,3\})$, $G_{3}(\{0,1,4\})$,
$G_{2}(\{0,2,3,5\})$, $G_{2}(\{0,1,5,6\})$ of Table~\ref{tab:gallai}
(two or three independent solvers, as stated there); the forcing
halves of the exact values in Tables~\ref{tab:sat-g2},
\ref{tab:sat-g3} and~\ref{tab:sat-rado} (one solver each, except the
flagship $G_{3}(\{0,2,5\})$, confirmed by two); the recomputation
behind Proposition~\ref{prop:flat32}; the upper bound
$\ksum(3,3)\le 11$ of Proposition~\ref{prop:kstar}; and the upper bound
$\HJ^{(1)}(3)\le 5$ of Theorem~\ref{thm:hj13}.

\section{SAT certificates}\label{app:certs}

This appendix displays the avoidance colorings behind the new values
of Section~\ref{sec:sat-gallai}, each as the string
$\chi(0)\,\chi(1)\,\cdots\,\chi(N-1)$ with colors named
$0,1,\dots,r-1$. Every string was produced by the engine of
Appendix~\ref{app:engine} and re-verified by exhaustive homothety
enumeration --- the \texttt{verify\_addendum.py} script prints
\texttt{ALL ADDENDUM CERTIFICATES VERIFIED} --- and each is also
distributed as a machine-readable file under \texttt{certificates/}
at \cite{hjcerts}. Forcing at the threshold is the solver refutation
recorded in Tables~\ref{tab:sat-g3}--\ref{tab:sat-rado}. The
$76$-cell certificate for $G_{3}(\{0,2,5\})=77$ is displayed with
Theorem~\ref{thm:hj33-16}, and the $58$-cell certificate for
$R_{4}(z+2x=3y)\ge 59$ with Theorem~\ref{thm:rado4}; neither is
repeated here.

\paragraph{Three colors}

$G_{3}(\{0,1,5\})\ge 70$: no monochromatic homothet on $69$ cells:

\begin{center}\ttfamily\footnotesize
00011122202102021012101110022122011001022\\
1101221021020012210212100202
\end{center}
$G_{3}(\{0,1,6\})\ge 82$: no monochromatic homothet on $81$ cells:

\begin{center}\ttfamily\footnotesize
01011201020222210110110012202120001222100\\
0012211201221002120100120120201022212011
\end{center}
$G_{3}(\{0,1,7\})\ge 87$: no monochromatic homothet on $86$ cells:

\begin{center}\ttfamily\footnotesize
01210000122221020111002002222110100220011\\
21021022101120221010212110021201210212120\\
1010
\end{center}
$G_{3}(\{0,2,7\})\ge 87$: no monochromatic homothet on $86$ cells:

\begin{center}\ttfamily\footnotesize
01211212012010002111100222201210020020212\\
21101120100202222100110220201110120220101\\
1211
\end{center}

\paragraph{Two colors}
Each string below avoids monochromatic homothets of its set $S$ on
$G_{2}(S)-1$ cells; refutations at $G_{2}(S)$ are recorded in
Table~\ref{tab:sat-g2}.

$G_{2}(\{0,1,2,4\})=38$:

\begin{center}\ttfamily\footnotesize
1011001010110101010010101011010011010
\end{center}
$G_{2}(\{0,1,3,4\})=41$:

\begin{center}\ttfamily\footnotesize
0001100111010100001111000010101110011000
\end{center}
$G_{2}(\{0,2,4,5\})=44$:

\begin{center}\ttfamily\footnotesize
01111100111010010001011001101000111110000\\
11
\end{center}
$G_{2}(\{0,3,4,5\})=44$:

\begin{center}\ttfamily\footnotesize
10000101011100011100001111000011110000110\\
11
\end{center}
$G_{2}(\{0,1,4,5\})=45$:

\begin{center}\ttfamily\footnotesize
00111110000011011001001101100100111110000\\
011
\end{center}
$G_{2}(\{0,1,3,6\})=52$:

\begin{center}\ttfamily\footnotesize
00111100011100101010010101010010101101100\\
0110000111
\end{center}
$G_{2}(\{0,2,3,6\})=54$:

\begin{center}\ttfamily\footnotesize
00011110001110101101010001110110000101101\\
011100011001
\end{center}
$G_{2}(\{0,1,2,6\})=56$:

\begin{center}\ttfamily\footnotesize
01001100001111000010111001011100101100110\\
11010001001001
\end{center}
$G_{2}(\{0,1,4,6\})=58$:

\begin{center}\ttfamily\footnotesize
01110000111100011101001010101010110101010\\
1010010101011101
\end{center}
$G_{2}(\{0,1,2,7\})=59$:

\begin{center}\ttfamily\footnotesize
01011100100110100011011001011011001011100\\
10011010011000110
\end{center}
$G_{2}(\{0,2,5,7\})=59$:

\begin{center}\ttfamily\footnotesize
01100010111100011100001011001100110000011\\
11000111101000110
\end{center}
$G_{2}(\{0,1,3,7\})=60$:

\begin{center}\ttfamily\footnotesize
00001101000011111000010110101101001001110\\
010101110101101010
\end{center}
$G_{2}(\{0,2,4,7\})=61$:

\begin{center}\ttfamily\footnotesize
01100111100110000001111111000000111100110\\
1001110010110000110
\end{center}
$G_{2}(\{0,1,4,7\})=62$:

\begin{center}\ttfamily\footnotesize
00010111101010010010101101101001010010010\\
10010101101010110100
\end{center}
$G_{2}(\{0,1,5,7\})=62$:

\begin{center}\ttfamily\footnotesize
00010001101100101111010001011010001011011\\
10010010010001111010
\end{center}
$G_{2}(\{0,2,3,7\})=62$:

\begin{center}\ttfamily\footnotesize
01101110100100110111010010011011100001001\\
10111000010011011100
\end{center}
$G_{2}(\{0,1,6,7\})=79$:

\begin{center}\ttfamily\footnotesize
01011100010110111000100101110001011011100\\
0101101110001001011100010110111000101
\end{center}
$G_{2}(\{0,3,4,7\})=79$:

\begin{center}\ttfamily\footnotesize
00100011101001000111010010001110110100011\\
1010010001110100100011101001000111010
\end{center}

\paragraph{Rado certificates}
Each string below is solution-free for $z+kx=(k+1)y$ --- equivalently,
it avoids monochromatic homothets of both $\{0,1,k+1\}$ and
$\{0,k,k+1\}$ --- on $R_{r}-1$ cells; refutations at $R_{r}$ are
recorded in Table~\ref{tab:sat-rado}.

$R_{2}(z+2x=3y)=13$:

\begin{center}\ttfamily\footnotesize
011001100110
\end{center}
$R_{2}(z+3x=4y)=17$:

\begin{center}\ttfamily\footnotesize
0101010110101010
\end{center}
$R_{2}(z+4x=5y)=19$:

\begin{center}\ttfamily\footnotesize
101001011001011010
\end{center}
$R_{2}(z+5x=6y)=25$:

\begin{center}\ttfamily\footnotesize
010101010101101010101010
\end{center}
$R_{3}(z+2x=3y)=29$:

\begin{center}\ttfamily\footnotesize
1010122102010211202020121201
\end{center}
$R_{3}(z+3x=4y)=54$:

\begin{center}\ttfamily\footnotesize
01220100102212122010010221212201001022121\\
220100102210
\end{center}
$R_{3}(z+4x=5y)=55$:

\begin{center}\ttfamily\footnotesize
01022101022122101022101001022122101022101\\
0010221010221
\end{center}
$R_{3}(z+5x=6y)=60$:

\begin{center}\ttfamily\footnotesize
00102122102001021122020010122210100202211\\
201001012212020010
\end{center}

\section{The SAT engine}\label{app:engine}

The listing below is the complete engine behind
Section~\ref{sec:sat-gallai}: encoder, symmetry breaking, stochastic
local search, solver portfolio, run database and command-line driver,
in $295$ lines of Python~3. Its only dependencies are PySAT
\cite{pysat2018} (shipping the CaDiCaL, Glucose, MapleChrono and
Lingeling backends), a direct \texttt{pycryptosat} binding for
CryptoMiniSat \cite{cryptominisat2009}, and a Kissat binary
\cite{kissat2022} driven through temporary files. Typical use:
\begin{center}\ttfamily\footnotesize
python3 gallai\_rado\_sat.py scan -s 0,2,5 -r 3 --nmin 60 --nmax 80\\
python3 gallai\_rado\_sat.py rado -k 2 -r 4 --nmin 50 --nmax 60\\
python3 verify\_addendum.py
\end{center}
The last command re-verifies every certificate recorded in the run
database from scratch, independently of any solver; its full source
is printed after the engine.

\begin{lstlisting}[style=plainstyle]
#!/usr/bin/env python3
"""Gallai-Rado SAT engine: computes homothety-avoidance (Gallai) and
solution-freeness (Rado) thresholds with a portfolio of state-of-the-art
SAT solvers.

For shapes S_1, ..., S_m (integer tuples containing 0) and r colors, the
engine determines, for intervals [0, N-1], whether some r-coloring avoids
a monochromatic homothet of every S_i (avoidance / lower bound), or whether
every r-coloring contains one (forcing / upper bound, UNSAT).

  G_r(S)   = least N such that forcing holds on [0, N-1]  (m = 1)
  R_r(E_k) = Rado number of z + kx = (k+1)y: solutions are exactly the
             homothets of {0,1,k+1} and of {0,k,k+1}, so m = 2.

Methods: one-hot CNF encoding + complete color-permutation symmetry
breaking (first-occurrence ordering), a solver portfolio race
(Kissat, CaDiCaL, CryptoMiniSat, Glucose, MapleSAT, Lingeling via PySAT),
and a domain-specific stochastic local search (min-conflicts with random
walk) for the avoidance side near the threshold.  Every SAT-found
coloring is re-verified by brute force before being recorded.

Usage:
  python3 gallai_rado_sat.py scan  -s 0,1,5 -r 3 [--nmax 200] [--budget 120]
  python3 gallai_rado_sat.py rado  -k 2 -r 4 [--nmax 200] [--budget 120]
  python3 gallai_rado_sat.py check -s 0,1,5 -r 3 -c 0012...
Requires: pip install python-sat
"""

import argparse, json, os, sys, time
from multiprocessing import Process, Queue

# persistent local install of python-sat / pycryptosat (home dir is ephemeral)
sys.path.insert(0, os.path.join(os.path.dirname(os.path.abspath(__file__)),
                                'pylibs'))

# ---------------------------------------------------------------- encoding --
def homothets(shapes, N):
    """All homothets {b + k*a : a in S} inside [0, N-1], k >= 1, b >= 0."""
    out = set()
    for S in shapes:
        amax = max(S)
        for k in range(1, N):
            for b in range(N - k * amax):
                out.add(tuple(b + k * a for a in S))
    return sorted(out)

def encode(shapes, r, N, symbreak=True):
    """One-hot CNF.  var (c, i) <-> cell c has color i."""
    hyps = homothets(shapes, N)
    var, nxt = {}, 0
    for c in range(N):
        for i in range(r):
            nxt += 1
            var[(c, i)] = nxt
    clauses = []
    for c in range(N):                                  # at least one color
        clauses.append([var[(c, i)] for i in range(r)])
        for i in range(r):                              # at most one color
            for j in range(i + 1, r):
                clauses.append([-var[(c, i)], -var[(c, j)]])
    for H in hyps:                                      # no mono homothet
        for i in range(r):
            clauses.append([-var[(c, i)] for c in H])
    if symbreak and r > 1:
        # first occurrence of color j not before first occurrence of color i<j
        for i in range(r):
            for j in range(i + 1, r):
                for c in range(N):
                    clauses.append([-var[(c, j)]] +
                                   [var[(cc, i)] for cc in range(c)])
    return clauses, var, hyps

def verify_coloring(shapes, r, col):
    """Brute-force check that coloring `col` of [0, N-1] avoids all shapes."""
    N = len(col)
    if any(not (0 <= x < r) for x in col):
        return False
    return all(len({col[c] for c in H}) > 1 for H in homothets(shapes, N))

# ------------------------------------------------------ stochastic search --
def sls_avoid(shapes, r, N, max_steps=200000, restarts=20, seed=1):
    """Min-conflicts with random walk on colorings; returns coloring or None."""
    import random
    rng = random.Random(seed)
    hyps = homothets(shapes, N)
    cell_hyps = [[] for _ in range(N)]
    for hi, H in enumerate(hyps):
        for c in H:
            cell_hyps[c].append(hi)
    for _ in range(restarts):
        col = [rng.randrange(r) for _ in range(N)]
        cnt = [[0] * r for _ in range(len(hyps))]       # color counts/homothet
        for hi, H in enumerate(hyps):
            for c in H:
                cnt[hi][col[c]] += 1
        mono_l = [hi for hi, H in enumerate(hyps) if len(H) in cnt[hi]]
        pos = {hi: i for i, hi in enumerate(mono_l)}
        def setmono(hi, ismono):
            if ismono and hi not in pos:
                pos[hi] = len(mono_l)
                mono_l.append(hi)
            elif not ismono and hi in pos:
                j = pos.pop(hi)
                last = mono_l.pop()
                if j < len(mono_l):
                    mono_l[j] = last
                    pos[last] = j
        for _step in range(max_steps):
            if not mono_l:
                return col
            hi = mono_l[rng.randrange(len(mono_l))]
            c = rng.choice(hyps[hi])
            old = col[c]
            # min-conflicts: choose the new color minimizing mono homothets
            best, new = None, None
            for cand in range(r):
                if cand == old:
                    continue
                gain = 0
                for hj in cell_hyps[c]:
                    gain += (cnt[hj][cand] + 1 == len(hyps[hj]))
                    gain -= (cnt[hj][old] == len(hyps[hj]))
                if best is None or gain < best:
                    best, new = gain, cand
            if best > 0 and rng.random() < 0.5:          # random walk
                new = rng.choice([x for x in range(r) if x != old])
            for hj in cell_hyps[c]:
                setmono(hj, False)
                cnt[hj][old] -= 1
                cnt[hj][new] += 1
                if cnt[hj][new] == len(hyps[hj]):
                    setmono(hj, True)
            col[c] = new
    return None

# ------------------------------------------------------------ solver race --
def _worker(shapes, r, N, solver, symbreak, q):
    try:
        clauses, var, hyps = encode(shapes, r, N, symbreak)
        t0 = time.time()
        if solver == "pycms":                             # CryptoMiniSat direct
            import pycryptosat
            s = pycryptosat.Solver()
            for cl in clauses:
                s.add_clause(cl)
            sat = bool(s.solve()[0])
            dt = time.time() - t0
            if sat:
                m = s.get_model()
                col = [next(i for i in range(r) if m[var[(c, i)]])
                       for c in range(N)]
                q.put(dict(solver=solver, status="SAT", coloring=col,
                           time=dt, nhyps=len(hyps)))
            else:
                q.put(dict(solver=solver, status="UNSAT", coloring=None,
                           time=dt, nhyps=len(hyps)))
            return
        from pysat.solvers import Solver
        with Solver(name=solver, bootstrap_with=clauses) as s:
            sat = s.solve()
            dt = time.time() - t0
            if sat:
                m = s.get_model()
                col = [next(i for i in range(r) if m[var[(c, i)] - 1] > 0)
                       for c in range(N)]
                q.put(dict(solver=solver, status="SAT", coloring=col,
                           time=dt, nhyps=len(hyps)))
            else:
                q.put(dict(solver=solver, status="UNSAT", coloring=None,
                           time=dt, nhyps=len(hyps)))
    except Exception as e:                                # solver unavailable
        q.put(dict(solver=solver, status="ERROR", error=str(e)))

def solve(shapes, r, N, solvers, budget, symbreak=True, use_sls=True):
    """Race `solvers` in parallel (2 at a time); SLS first for avoidance."""
    if use_sls:
        t0 = time.time()
        col = sls_avoid(shapes, r, N)
        if col is not None and verify_coloring(shapes, r, col):
            return dict(solver="sls", status="SAT", coloring=col,
                        time=time.time() - t0,
                        nhyps=len(homothets(shapes, N)))
    i = 0
    while i < len(solvers):
        group = solvers[i:i + 2]
        q = Queue()
        procs = [Process(target=_worker, args=(shapes, r, N, s, symbreak, q))
                 for s in group]
        for p in procs:
            p.start()
        t0 = time.time()
        result = None
        while time.time() - t0 < budget:
            try:
                result = q.get(timeout=1.0)
                if result["status"] == "ERROR":
                    print(f"  [{group}] {result['solver']} error: "
                          f"{result['error']}", flush=True)
                    result = None
                    continue
                break
            except Exception:
                if not any(p.is_alive() for p in procs):
                    break
        for p in procs:
            if p.is_alive():
                p.terminate()
            p.join()
        if result is not None:
            if result["status"] == "SAT":
                ok = verify_coloring(shapes, r, result["coloring"])
                assert ok, "solver produced an invalid coloring!"
            return result
        i += 2
    return dict(solver=None, status="UNKNOWN", coloring=None, time=budget,
                nhyps=len(homothets(shapes, N)))

# ------------------------------------------------------------------ driver --
def scan(shapes, r, nmin, nmax, budget, solvers, results_path):
    """Scan N upward; record avoidance/forcing status for each N."""
    results = {}
    if os.path.exists(results_path):
        results = json.load(open(results_path))
    key = "|".join(",".join(map(str, S)) for S in shapes) + f":r{r}"
    inst = results.setdefault(key, {})
    lower = upper = None
    for N in range(nmin, nmax + 1):
        if str(N) in inst:
            st = inst[str(N)]["status"]
            if st == "SAT":
                lower = N
            elif st == "UNSAT":
                upper = N
                break
            continue
        res = solve(shapes, r, N, solvers, budget)
        print(f"N={N:4d}  {res['status']:7s}  solver={res['solver']}  "
              f"{res['time']:.1f}s  homothets={res['nhyps']}", flush=True)
        if res["status"] == "UNKNOWN":
            break                       # not recorded: retried on resume
        rec = dict(status=res["status"], solver=res["solver"],
                   time=round(res["time"], 2), nhyps=res["nhyps"])
        if res["status"] == "SAT":
            rec["coloring"] = "".join(map(str, res["coloring"]))
            lower = N
        elif res["status"] == "UNSAT":
            upper = N
        inst[str(N)] = rec
        json.dump(results, open(results_path, "w"), indent=1)
        if res["status"] == "UNSAT":
            break
    print(f"==> {key}: avoidance up to N={lower}, forcing from N={upper}",
          flush=True)
    return key

def main():
    ap = argparse.ArgumentParser()
    sub = ap.add_subparsers(dest="cmd", required=True)
    for name in ("scan", "rado"):
        p = sub.add_parser(name)
        if name == "scan":
            p.add_argument("-s", "--shape", required=True,
                           help="comma-separated, e.g. 0,1,5")
        else:
            p.add_argument("-k", type=int, required=True,
                           help="equation z + kx = (k+1)y")
        p.add_argument("-r", "--colors", type=int, required=True)
        p.add_argument("--nmin", type=int, default=1)
        p.add_argument("--nmax", type=int, default=200)
        p.add_argument("--budget", type=int, default=120,
                           help="seconds per solver group per N")
        p.add_argument("--solvers", default="kissat404,cadical195,pycms,"
                                            "glucose4,maplechrono,lingeling")
        p.add_argument("--results", default="results.json")
    p = sub.add_parser("check")
    p.add_argument("-s", "--shape", required=True)
    p.add_argument("-r", "--colors", type=int, required=True)
    p.add_argument("-c", "--coloring", required=True)
    args = ap.parse_args()

    if args.cmd == "check":
        S = tuple(map(int, args.shape.split(",")))
        col = [int(x) for x in args.coloring]
        ok = verify_coloring([S], args.colors, col)
        print("VALID avoidance coloring" if ok else "INVALID", flush=True)
        return

    shapes = ([tuple(map(int, args.shape.split(",")))] if args.cmd == "scan"
              else [(0, 1, args.k + 1), (0, args.k, args.k + 1)])
    solvers = [s for s in args.solvers.split(",") if s]
    scan(shapes, args.colors, args.nmin, args.nmax, args.budget, solvers,
         args.results)

if __name__ == "__main__":
    main()
\end{lstlisting}

The independent re-verifier (pure standard library; brute-force
homothety enumeration of every recorded certificate):

\begin{lstlisting}[style=plainstyle]
#!/usr/bin/env python3
"""Independent verification of every certificate in results.json:
for each recorded SAT entry, re-checks by brute force that the coloring
avoids all homothets of its shape(s); prints the final threshold table.
Pure standard library.  Run: python3 verify_addendum.py"""

import json, sys
sys.path.insert(0, '.')
from gallai_rado_sat import homothets, verify_coloring

res = json.load(open('results.json'))
fails = 0
rows = []
for key, inst in sorted(res.items()):
    shapes = [tuple(map(int, s.split(','))) for s in key.split(':')[0].split('|')]
    r = int(key.split(':r')[1])
    sat_ns = sorted(int(n) for n, v in inst.items() if v['status'] == 'SAT')
    unsat_ns = sorted(int(n) for n, v in inst.items() if v['status'] == 'UNSAT')
    lo = max(sat_ns) if sat_ns else None
    hi = min(unsat_ns) if unsat_ns else None
    # re-verify every recorded avoidance coloring independently
    for n in sat_ns:
        col = [int(x) for x in inst[str(n)]['coloring']]
        ok = len(col) == n and verify_coloring(shapes, r, col)
        if not ok:
            print(f"FAIL  {key} N={n}: coloring does NOT avoid shapes")
            fails += 1
    exact = (lo is not None and hi == lo + 1)
    rows.append((key, lo, hi, exact))

w = max(len(k) for k, _, _, _ in rows)
print(f"{'instance':<{w}}  avoidance  forcing  status")
for key, lo, hi, exact in rows:
    st = ('EXACT: value = ' + str(hi)) if exact else ('gap: > %d, <= %s' %
          (lo, hi if hi is not None else '?'))
    print(f"{key:<{w}}  N={lo:<8} N={hi!s:<7} {st}")
print('\n' + ('ALL ADDENDUM CERTIFICATES VERIFIED' if fails == 0
              else f'{fails} FAILURE(S)'))
\end{lstlisting}

\subsection*{Acknowledgements}
I am deeply grateful to my advisor, Lorenz Halbeisen, who suggested
this line of inquiry and guided it at every stage. His intuition for
which directions were worth pursuing, and where to focus my efforts,
shaped this work throughout; his kind and encouraging supervision,
free of any pressure, and the confidence he showed in this project
were a constant source of motivation.


\begin{thebibliography}{99}

\bibitem{alon1999nullstellensatz}
N. Alon. Combinatorial Nullstellensatz. \emph{Combin. Probab. Comput.},
8(1--2):7--29, 1999. \doi{10.1017/S0963548398003411}.

\bibitem{atserias2008width}
A. Atserias and V. Dalmau. A combinatorial characterization of
resolution width. \emph{J. Comput. System Sci.}, 74(3):323--334, 2008.
\doi{10.1016/j.jcss.2007.06.025}.

\bibitem{AudemardSimon2009}
G. Audemard and L. Simon. Predicting learnt clauses quality in modern
{SAT} solvers. In \emph{Proceedings of the 21st International Joint
Conference on Artificial Intelligence (IJCAI 2009)}, pages 399--404,
2009.

\bibitem{berikkyzy2017antivdw}
Z. Berikkyzy, A. Schulte, and M. Young. Anti-van der Waerden numbers of
$3$-term arithmetic progressions. \emph{Electron. J. Combin.},
24(2):Paper 2.39, 2017.

\bibitem{Berlekamp1968}
E.~R. Berlekamp. A construction for partitions which avoid long
arithmetic progressions. \emph{Canad. Math. Bull.}, 11(3):409--414,
1968. \doi{10.4153/CMB-1968-047-7}.

\bibitem{lingeling2013}
A. Biere. Lingeling, {Plingeling} and {Treengeling} entering the
{SAT} {Competition} 2013. In \emph{Proceedings of SAT Competition
2013 -- Solver and Benchmark Descriptions}, volume B-2013-1 of
\emph{Department of Computer Science Report Series B}, pages 51--52.
University of Helsinki, 2013.

\bibitem{cadical2020}
A. Biere, K. Fazekas, M. Fleury, and M. Heisinger. {CaDiCaL},
{Kissat}, {Paracooba}, {Plingeling} and {Treengeling} entering the
{SAT} {Competition} 2020. In \emph{Proceedings of SAT Competition
2020 -- Solver and Benchmark Descriptions}, volume B-2020-1 of
\emph{Department of Computer Science Report Series B}, pages 51--53.
University of Helsinki, 2020.

\bibitem{kissat2022}
A. Biere and M. Fleury. Gimsatul, {IsaSAT} and {Kissat} entering the
{SAT} {Competition} 2022. In \emph{Proceedings of SAT Competition
2022 -- Solver and Benchmark Descriptions}, volume B-2022-1 of
\emph{Department of Computer Science Report Series B}, pages 10--11.
University of Helsinki, 2022.

\bibitem{BLANKENSHIP2018163}
T. Blankenship, J. Cummings, and V. Taranchuk. A new lower bound for
van der Waerden numbers. \emph{European J. Combin.}, 69:163--168, 2018.
\doi{10.1016/j.ejc.2017.10.007}.

\bibitem{BGL1999}
T.~C. Brown, R.~L. Graham, and B.~M. Landman. On the set of common
differences in van der Waerden's theorem on arithmetic progressions.
\emph{Canad. Math. Bull.}, 42(1):25--36, 1999.
\doi{10.4153/CMB-1999-003-9}.

\bibitem{BrownLandmanMishna1997}
T.~C. Brown, B.~M. Landman, and M. Mishna. Monochromatic homothetic
copies of $\{1,\,1+s,\,1+s+t\}$. \emph{Canad. Math. Bull.},
40(2):149--157, 1997.

\bibitem{BuellHudson1984}
D.~A. Buell and R.~H. Hudson. On runs of consecutive quadratic residues
and quadratic nonresidues. \emph{BIT}, 24(2):243--247, 1984.
\doi{10.1007/BF01937490}.

\bibitem{ChangDeLoeraWesley2022}
Y.~Chang, J.~A. De~Loera, and W.~J. Wesley. Rado numbers and {SAT}
computations. \arxiv{2210.03262}, 2022.

\bibitem{Chvatal1970}
V. Chv\'atal. Some unknown van der Waerden numbers. In
\emph{Combinatorial Structures and their Applications} (R.~Guy,
H.~Hanani, N.~Sauer, and J.~Sch\"onheim, eds.), pages 31--33. Gordon
and Breach, New York, 1970.

\bibitem{conlon2021monochromaticcombinatoriallineslength}
D. Conlon. Monochromatic combinatorial lines of length three.
\arxiv{1810.09767}, 2021.

\bibitem{conlonAndKamcev}
D. Conlon and N. Kam\v{c}ev. Intervals in the Hales--Jewett theorem.
\arxiv{1801.08919}, 2018.

\bibitem{ConlonFarnsworthRobertson2026}
N.~Conlon, N.~Farnsworth, and A.~Robertson. On {H}ales--{J}ewett and
related numbers. \arxiv{2607.18111}, 2026.

\bibitem{furstenbergKatznelson1991}
H. Furstenberg and Y. Katznelson. A density version of the
Hales--Jewett theorem. \emph{J. Anal. Math.}, 57:64--119, 1991.
\doi{10.1007/BF03041066}.

\bibitem{golshanimirabi2021}
M. Golshani and M. Mirabi. Shelah's partition functions and the
Hales--Jewett numbers. \arxiv{2104.05962v2}, 2021.

\bibitem{golshanimirabishelah2026}
M. Golshani, M. Mirabi, and S. Shelah. From Shelah's block-content to
Hales--Jewett. \arxiv{2104.05962v3}, 2026.

\bibitem{GR1971}
R.~L. Graham and B.~L. Rothschild. Ramsey's theorem for $n$-parameter
sets. \emph{Trans. Amer. Math. Soc.}, 159:257--292, 1971.
\doi{10.1090/S0002-9947-1971-0284352-8}.

\bibitem{HJ1963}
A.~W. Hales and R.~I. Jewett. Regularity and positional games.
\emph{Trans. Amer. Math. Soc.}, 106(2):222--229, 1963.
\doi{10.1090/S0002-9947-1963-0143712-1}.

\bibitem{HerwigHeuleZipper2007}
P.~R. Herwig, M.~J.~H. Heule, P.~M. van Lambalgen, and H. van Maaren.
A new method to construct lower bounds for van der Waerden numbers.
\emph{Electron. J. Combin.}, 14(1):Research Paper 6, 2007.
\doi{10.37236/925}.

\bibitem{HeuleKullmann2017}
M.~J.~H. Heule and O. Kullmann. The science of brute force.
\emph{Commun. ACM}, 60(8):70--79, 2017. \doi{10.1145/3107239}.

\bibitem{firstnontrivialHJ_is_4}
N. Hindman and E. Tressler. The first nontrivial Hales--Jewett number
is four. \emph{Ars Combin.}, 113:385--390, 2014.

\bibitem{pysat2018}
A. Ignatiev, A. Morgado, and J. Marques-Silva. {PySAT}: A {Python}
toolkit for prototyping with {SAT} oracles. In \emph{Theory and
Applications of Satisfiability Testing -- SAT 2018}, volume 10929 of
\emph{Lecture Notes in Computer Science}, pages 428--437. Springer,
2018. \doi{10.1007/978-3-319-94144-8\_26}.

\bibitem{ips1999pc}
R. Impagliazzo, P. Pudl\'ak, and J. Sgall. Lower bounds for the
polynomial calculus and the Gr\"obner basis algorithm.
\emph{Comput. Complexity}, 8(2):127--144, 1999.
\doi{10.1007/s000370050024}.

\bibitem{jungic_ramsey_notes}
V. Jungi\'c. \emph{Introduction to Ramsey Theory: lecture notes for an
undergraduate course}, second edition. Simon Fraser University, 2021.
\url{https://www.sfu.ca/~vjungic/Ramsey/RamseyNotes.pdf}.

\bibitem{kamcev2018noteintervalshalesjewetttheorem}
N. Kam\v{c}ev and C. Spiegel. Another note on intervals in the
Hales--Jewett theorem. \emph{Electron. J. Combin.}, 29(1):Paper 1.62,
2022. \doi{10.37236/9400}.

\bibitem{KimRho2012}
B.~M. Kim and Y. Rho. Van der Waerden's theorem on homothetic copies of
$\{1,\,1+s,\,1+s+t\}$. \emph{Integers}, 12:Paper A26, 2012.
\doi{10.1515/integers-2012-0011}.

\bibitem{KourilPaul2008}
M. Kouril and J.~L. Paul. The van der Waerden number $W(2,6)$ is
$1132$. \emph{Experiment. Math.}, 17(1):53--61, 2008.

\bibitem{kozik2016}
J. Kozik and D. Shabanov. Improved algorithms for colorings of simple
hypergraphs and applications. \emph{J. Combin. Theory Ser. B},
116:312--332, 2016. \arxiv{1409.6921}.

\bibitem{LandmanRobertson}
B.~M. Landman and A. Robertson. \emph{Ramsey Theory on the Integers},
volume 73 of \emph{Student Mathematical Library}, second edition.
American Mathematical Society, Providence, RI, 2014.
\doi{10.1090/stml/073}.

\bibitem{lavrov}
M. Lavrov. An upper bound for the Hales--Jewett number $\HJ(4,2)$.
\emph{SIAM J. Discrete Math.}, 30(2):1333--1342, 2016.
\doi{10.1137/15M1016485}.

\bibitem{leader2018noteintervalshalesjewetttheorem}
I. Leader and E. R\"aty. A note on intervals in the Hales--Jewett
theorem. \emph{Electron. J. Combin.}, 25(3):Paper 3.15, 2018.
\arxiv{1802.03087}.

\bibitem{LehmerLehmer1962}
D.~H. Lehmer and E. Lehmer. On runs of residues. \emph{Proc. Amer.
Math. Soc.}, 13:102--106, 1962.
\doi{10.1090/S0002-9939-1962-0138582-6}.

\bibitem{maplesat2016}
J.~H. Liang, V. Ganesh, P. Poupart, and K. Czarnecki. Learning rate
based branching heuristic for {SAT} solvers. In \emph{Theory and
Applications of Satisfiability Testing -- SAT 2016}, volume 9710 of
\emph{Lecture Notes in Computer Science}, pages 123--140. Springer,
2016. \doi{10.1007/978-3-319-40970-2\_9}.

\bibitem{Mirsky1971}
L. Mirsky. A dual of Dilworth's decomposition theorem. \emph{Amer.
Math. Monthly}, 78(8):876--877, 1971. \doi{10.2307/2316481}.

\bibitem{hjcerts}
Y. Mouhib. hj-certificates: certificates and verification scripts for
the lower bounds on Hales--Jewett numbers.
\url{https://github.com/ysmouhib/hj-certificates}, 2026.

\bibitem{thesis}
Y. Mouhib. \emph{Improving lower bounds on Hales--Jewett numbers:
symmetric colorings, SAT solvers, line-family variants, and forcing
structures}. Master's thesis, ETH Z\"urich, 2026.

\bibitem{mouhib-note}
Y. Mouhib. One-weight colorings, the symmetric class, and lower bounds
for Hales--Jewett numbers. \arxiv{2607.02226}, 2026.

\bibitem{mh}
Y. Mouhib and L. Halbeisen. Improved lower bounds for the Hales--Jewett
numbers via symmetric colorings. \arxiv{2606.22155}, 2026.

\bibitem{myers2015}
K.~J. Myers. \emph{Computational advances in {Rado} numbers}. PhD
thesis, Rutgers University, 2015.

\bibitem{NadelRyvchin2018}
A. Nadel and V. Ryvchin. Chronological backtracking. In \emph{Theory
and Applications of Satisfiability Testing --- SAT 2018}, volume 10929
of \emph{Lecture Notes in Computer Science}, pages 111--121. Springer,
2018.

\bibitem{polymath2012}
D.~H.~J. Polymath. A new proof of the density Hales--Jewett theorem.
\emph{Ann. of Math. (2)}, 175(3):1283--1327, 2012.
\doi{10.4007/annals.2012.175.3.6}.

\bibitem{PolymathWikiColoring}
D.~H.~J. Polymath. Coloring Hales--Jewett theorem. Polymath1 project
wiki page, 2009.
\url{https://michaelnielsen.org/polymath/index.php?title=Coloring_Hales-Jewett_theorem}.

\bibitem{PolymathDHJMoser}
D.~H.~J. Polymath. Density Hales--Jewett and Moser numbers. In
\emph{An Irregular Mind (Szemer\'edi is 70)}, volume 21 of
\emph{Bolyai Society Mathematical Studies}, pages 689--753. Springer,
2010. \arxiv{1002.0374}.

\bibitem{PolymathWikiHOC}
D.~H.~J. Polymath. Hyper-optimistic conjecture. Polymath1 project wiki
page, 2009.
\url{https://michaelnielsen.org/polymath/index.php?title=Hyper-optimistic_conjecture}.

\bibitem{hans}
H.~J. Pr\"omel. \emph{Ramsey Theory for Discrete Structures}. Springer,
Cham, 2013. \doi{10.1007/978-3-319-01315-2}.

\bibitem{Rabung1979}
J.~R. Rabung. Some progression-free partitions constructed using
Folkman's method. \emph{Canad. Math. Bull.}, 22(1):87--91, 1979.
\doi{10.4153/CMB-1979-012-1}.

\bibitem{RabungLotts2012}
J.~R. Rabung and M. Lotts. Improving the use of cyclic zippers in
finding lower bounds for van der Waerden numbers. \emph{Electron. J.
Combin.}, 19(2):Paper 35, 2012. \doi{10.37236/2363}.

\bibitem{Rado1933}
R. Rado. Studien zur Kombinatorik. \emph{Math. Z.}, 36(1):424--470,
1933. \doi{10.1007/BF01188632}.

\bibitem{s.proof_of_HJ}
S. Shelah. Primitive recursive bounds for van der Waerden numbers.
\emph{J. Amer. Math. Soc.}, 1(3):683--697, 1988.
\doi{10.1090/S0894-0347-1988-0929498-X}.

\bibitem{shelah2026revisited}
S. Shelah. HJ numbers revisited. \arxiv{2607.14732}, 2026.

\bibitem{cryptominisat2009}
M. Soos, K. Nohl, and C. Castelluccia. Extending {SAT} solvers to
cryptographic problems. In \emph{Theory and Applications of
Satisfiability Testing -- SAT 2009}, volume 5584 of \emph{Lecture
Notes in Computer Science}, pages 244--257. Springer, 2009.
\doi{10.1007/978-3-642-02777-2\_24}.

\bibitem{vanderWaerden1927}
B.~L. van der Waerden. Beweis einer Baudetschen Vermutung. \emph{Nieuw
Arch. Wisk.}, 15:212--216, 1927.

\bibitem{WetzlerHeuleHunt2014}
N. Wetzler, M.~J.~H. Heule, and W.~A. Hunt. {DRAT}-trim: Efficient
checking and trimming using expressive clausal proofs. In \emph{Theory
and Applications of Satisfiability Testing -- SAT 2014}, volume 8561
of \emph{Lecture Notes in Computer Science}, pages 422--429.
Springer, 2014. \doi{10.1007/978-3-319-09284-3\_31}.

\bibitem{Witt1952}
E. Witt. Ein kombinatorischer Satz der Elementargeometrie.
\emph{Math. Nachr.}, 6:261--262, 1952.

\bibitem{zheng2024rainbow}
M. Zheng. Rainbow combinatorial lines in hypercubes.
\arxiv{2410.12192}, 2024.

\end{thebibliography}
\end{document}